\documentclass[11pt,reqno]{amsart}

\usepackage{fullpage,cite}
\usepackage{amsfonts}
\usepackage{amssymb}
\usepackage{times}
\usepackage{algorithm}
\usepackage{algorithmic}
\usepackage[utf8]{inputenc}
\usepackage[portuges,english]{babel} 
\usepackage{hhline}
\usepackage{url}
\usepackage{graphicx}
\usepackage{tabularx}
\usepackage{makecell}
\usepackage{fancyvrb}
\usepackage[font=tiny,labelfont=bf]{caption}
\usepackage[foot]{amsaddr}

\theoremstyle{remark}
\usepackage[toc,page]{appendix}

\title{CREAM: a Package to Compute [Auto, Endo, Iso, Mono, Epi]-morphisms, Congruences, Divisors and More for Algebras of Type $(2^n,1^n)$}
\author{Jo\~ao Ara\'ujo}
\address{Mathematics Department, Faculty of Science and Technology, Universidade Nova de Lisboa,
 Portugal}
\email{jj.araujo@fct.unl.pt}
\author{Rui Barradas Pereira}
\address{Department of Science and Technology, Universidade Aberta, Lisbon, Portugal }
\email{rmbper@gmail.com}
\author{Wolfram Bentz}
\address{Department of Science and Technology, Universidade Aberta, Lisbon, Portugal }
\email{wolfram.bentz@uab.pt}
\author{Choiwah Chow}
\address{Department of Science and Technology, Universidade Aberta, Lisbon, Portugal }
\email{choiwah.chow@gmail.com}
\author{Jo\~ao Ramires}
\email{joao.j.ramires@gmail.com}
\address{Department of Science and Technology, Universidade Aberta, Lisbon, Portugal }
\author{Luis Sequeira}
\address{Mathematics Department, Universidade de Lisboa, Lisbon, Portugal}
\email{lfsequeira@fc.ul.pt}
\author{Carlos Sousa}
\address{Department of Science and Technology, Universidade Aberta, Lisbon, Portugal }
\email{cfmsousa@sapo.pt}

\begin{document}

\maketitle

\begin{abstract}
The CREAM GAP package computes automorphisms, congruences,  endomorphisms and subalgebras 
of algebras with an arbitrary number of binary and unary operations; it also decides if between two such algebras there exists a monomorphism, an epimorphism, an isomorphism or if one is a divisor of the other. Thus it finds those objects for almost all algebras used in practice (groups, quasigroups in their various signatures, semigroups possibly with many unary operations, fields, semi-rings, quandles, logic algebras, etc). 

As a one-size-fits-all package, it only relies on universal algebra theorems, without taking advantage of specific theorems about, eg, groups or semigroups to reduce  the search space.  Canon and Holt produced very fast code to compute automorphisms of groups that outperform CREAM on orders larger than 128. Similarly, Mitchell et al. take advantage of deep theorems to compute automorphisms and congruences of completely 0-simple semigroups in a very efficient manner. However these domains  (groups of order above 128 and completely 0-simple semigroups) are among the very few examples of GAP code faster than our general purpose package CREAM. For the overwhelming majority of other classes of algebras, either ours is the first code  computing the above mentioned objects, or the existing algorithms are outperformed by CREAM, in some cases by several orders of magnitude. 

To get this performance, CREAM uses a mixture of universal algebra algorithms together with GAP coupled with artificial intelligence theorem proving tools (AITP)   and very delicate C implementations. As an example of the latter, we re-implement Freese's very clever algorithm for computing congruences in universal algebras, in a way that outperforms all other known implementations.  
\end{abstract}

\section{Introduction}
\label{intro section}

Investigation of automorphism groups of mathematical structures is one
of the classical algebraic problems. A cornerstone was the work of Evariste Galois, but its impact goes far beyond. In the words of P. J. Cameron \cite{Cameron1983}:
\begin{quote}
\emph{
In the famous Erlanger Programme, Klein proposed that geometry is the study of symmetry;
more precisely, the geometric properties of an object are those which are invariant under all automorphisms of the objects. ($\ldots$)
According to Artin, ``the investigation of symmetries of a given mathematical structure has always yielded the most powerful results.''}
\end{quote}
These ideas, expressed by Klein in 1872 \cite{Klein1872} and by Artin in 1957  \cite{Artin1957},
give some insight on why
 the computation  of automorphisms of various mathematical structures has been such an attractive topic for so many decades.  It should be observed here that this effort is closely linked to an
omnipresent  problem in mathematics -- stated explicitly by Ulam \cite{ulam} --  of describing the mathematical objects whose endomorphisms encode {\em all the relevant information} about them:  

$$\mbox{End}(A)\cong \mbox{End}(B)\Rightarrow A\cong B.$$ 

The connection between endomorphisms and congruences is well known, but the importance of congruences goes far beyond a tool to get endomorphisms; in fact their importance in modern algebra can hardly be overestimated.  

Congruences have close ties to the structure of algebras - for example, the subdirectly irreducible algebras are precisely the ones that possess a unique minimal nontrivial congruence; direct irreducibility is also characterized by properties of congruences. 
In Universal Algebra, a variety is a class of algebras of the same type, that satisfy a certain set of identities, or equivalently, by a celebrated theorem of Garrett Birkhoff (see \cite[Theorem II.11.9]{BURRIS}), a class closed under taking subalgebras, homomorphic images and direct products. Each variety is also generated by its subdirectly irreducible algebras.
Properties of congruence lattices of algebras over a variety are closely related to the identities satisfied in that variety. 
This gave rise to what is usually referred to the theory of Maltsev conditions - see, for example, \cite{GT84}.
The commutator of normal subgroups was generalized for algebras, first those in a congruence-permutable variety \cite{Smith76}, then for any algebra in a congruence modular variety \cite{HH79}; see \cite{FMcK87}. Thus notions like Abelianness or nilpotency can be applied meaningfully in a much wider context, and many fruitful consequences have been derived, for example in the context of natural dualities \cite{CD98}. 
Tame Congruence Theory \cite{HMcK88} emerged as a powerful tool to study the structure of finite algebras, and continues to provide a big leverage in  the study of locally finite varieties. These examples point to the critical importance of the study of congruences in algebra.

As observed above,  congruences are closely related to endomorphisms and the latter form a monoid (semigroup with identity) that very often encodes very relevant information about the original object and can be investigated with  the increasingly powerful techniques developed by experts in semigroup theory (especially since they started to massively apply the classification of finite simple groups \cite{JA1,JA2,JA3,JA4,JA5,JA6,JA6b,JA7,JA8,JA9,JA10,JA11,JA12,JA13,JA14}). 

Everything said above should make the case in favour of a general centralized and very effective GAP tool that computes congruences, automorphisms and endomorphisms of general algebras of type $(2^m,1^n)$. We chose these since the overwhelming majority of algebras used in practice have this type: it is a very rare occurrence that an algebraic structure actually uses ternary or higher arity operations (see \cite{BURRIS}).

The problem is that computing [endo]automorphisms/congruences of general algebras is a difficult task and hence
the majority of existing tools target specific classes  (groups, semigroups, quasigroups, etc.) in order to take advantage of the domain's theorems. Probably the best example is the GAP code to find automorphisms of groups, devised by  J. J. Cannon and D. Holt   \cite{cannon-autos-2003}, that takes advantage of many deep theorems and builds upon many past computational optimizations to quickly get auxiliary objects.  
In contrast, since we want our algorithms to apply for general algebras of type  $(2^m,1^n)$, we cannot rely on theorems that only hold in very specific cases. The challenge of this project  is thus to provide an effective tool to compute the objects  for finite algebras of type $(2^m,1^n)$ that works for all algebras of this type and does not fall far behind the more dedicated tools that take advantage of very specific domain theorems. 

The final result is that  CREAM compares favourably even with the specialized Cannon and Holt tool referred to above: for groups up to order 128, our general purpose tool finds the automorphisms of a group faster than their very optimized code. 
For most other classes of algebras, our tool outperforms specialized code, in some cases by large margins. This was achieved by a combination of results and ideas from different parts of computational algebra,  artificial intelligence theorem proving, and optimized C implementations.  Regarding congruences, we take advantage of a new implementation of  Freese's algorithm \cite{FREESE}, an algorithm that stands out by its quality and generality; for automorphisms, after many different tests and approaches, the most effective way turned out to be a generalization of the ideas developed in \cite{LOOPS3_4_1}. Once possessing a very good tool to compute automorphisms, finding [mono,epi,iso]-morphisms is straightforward. 

Finally, our tool also handles divisors. Given two algebras, $A$ and $B$, we say that $B$ divides $A$ if $B$ is a homomorphic image of a subalgebra of $A$. In different words, $B$ divides $A$ if  there is a subalgebra $C$ of $A$ and a congruence $\rho$ of $C$ such that $C/\rho$ is isomorphic to $B$. Like congruences, divisors are very important tools. For example, the celebrated Krohn--Rhodes theorem states that every finite semigroup $S$ divides a wreath product of finite simple groups, each dividing
$S$, together with copies of the $3$-element monoid $\{1,a,b\}$, where $ba=aa=a$, $ab=bb=b$ (for details see \cite{steinberg}). The key difficulty when finding divisors is to find the subalgebras of a given algebra; CREAM has an algorithm to compute them.

There are
numerous papers concerning the automorphism groups of particular
classes of algebras, for example, Schreier \cite{Sc36} and Mal'cev
\cite{Ma52} described all automorphisms of the semigroup of all mappings from a set to itself. Similar results have been obtained for various other structures
such as orders, equivalence relations, graphs, and hypergraphs; see
the survey papers \cite{Mo83} and \cite{Mo01}. More examples are
provided, among others, by Glusk\v{\i}n \cite{Gl65}, 
 Ara\'ujo and Konieczny \cite{ArKo03}, \cite{joao1}, and \cite{joao2}, Fitzpatrick and
Symons \cite{fitz}, Levi \cite{Le85} and  \cite{Le87}, Liber
\cite{Li53}, Magill \cite{Ma67}, Schein \cite{Sc70}, Sullivan
\cite{Su75}, and \v{S}utov \cite{Su61}.

Fast algorithms exist to find the congruences and automorphisms of groups and semigroups. For example, with GAP \cite{GAP4}, 
the package SEMIGROUPS includes the function CongruencesOfSemigroup that returns the congruences of a given semigroup. This function only works for a limited set of semigroups, belonging to the classes Simple, Brandt, Group, Zero Simple or Rectangular Band, and only for Rees Matrix semigroups and Rees Zero Matrix semigroups is it highly efficient. On a different platform, UACalc \cite{UACALC} supports a wider range of algebras. UACalc also implements the Freeze algorithm that is the basis for the calculation of congruences described here. 

Regarding automorphisms in GAP, many special properties of groups are used to implement the function {\em AutomorphismGroup} to efficiently find the automorphism group of a given group. However, this specialized method, while extremely efficient, works only on groups and not on any other more general algebraic structures, such as quasigroups, semigroups and magmas. Likewise, the Loops package \cite{LOOPS3_4_1} in GAP provides another version of the function {\em AutomorphismGroup} to compute all the automorphisms of a given quasigroup.  

To date, there are no known implementations of general functions for finding divisors, congruences and  [endo, auto, epi, mono, iso]-morphisms of magmas or algebras of type $(2^m, 1^n)$ in GAP. To fill the void, the CREAM package implements such functions in GAP with part of the code written in C for performance. 

This article is composed of 6 sections: this introduction, a short description of the mathematical concepts and background, the description of the algorithms used in the CREAM package, a comparative discussion of the performance of the package, a list of applications of the package and finally a short conclusion.

\section{Mathematical Background}
\label{math section}

\subsection{Algebra of Type ($2^m$,$1^n$)}
\label{subsection:algebra_type_2}

An Algebra in the sense of Universal Algebra is an algebraic structure consisting of a set $A$ together with a collection of operations on $A$.
(If we want to be more precise,  an algebra also contains an indexed scheme that references the operations, but we are not going to enter those technical details;
for a complete definition see \cite{BURRIS}.)

An $n$-ary operation on $A$ is a function that takes an $n$-tuple of elements from $A$ and returns a single element of $A$, that is, a function from $A^{n}$ to $A$. The number $n$ is called the \emph{arity} of the operation. For the scope of this package we are only considering operations with arity $1$ or $2$, i.e., unary and binary operations. An algebra of type ($2^m$,$1^n$) is a universal algebra with $m$ binary and $n$ unary operations. 

The package represents a finite Universal Algebra as a list of operations. An operation of arity $1$ is represented by a vector, while an operation of arity $2$ is represented by a square matrix. The underlying set $A$ is implicitly specified by the vector or matrix sizes, which need to agree for a valid representation.  If this size is $d$, the algebra is defined on the set $A=\{1,\dots,d\}$. We can safely ignore the ambiguity this introduces in the (uninteresting) case of an algebra without operations. Each vector or matrix describes the corresponding operation by listing all images in the obvious way. The following is an example of a representation for an algebra with a unary and a binary operation. 

\begin{verbatim}
[ [3, 1, 2], [ [1, 2, 3], [2, 3, 1], [3, 1, 2] ] ]
\end{verbatim}

For those CREAM functions that involve more than one algebra, the operations of the algebras need to be aligned by a common index scheme. In our representation, this index is indirectly provided by the order in which the operations are listed. The involved algebras need to be compatible, that is, they have the same number of operations of each arity, and the operations  are listed in the same position in the algebra representation. 

\subsection{Partition}
\label{subsection:partition}

A partition of a set $A$ is a collection of non-empty subsets of $A$, such that every element of $A$ is included in exactly one subset \cite{BURRIS}. Each subset in a partition is called a block or part. For example, $\{\{1,3\},\{2,4\}\}$ is a 2-blocks partition of the set $\{1,2,3,4\}$. 

A block of a partition is efficiently represented as a tree. A partition is represented by a forest, that is, a disjoint union of trees. The forest is physically presented in the computer memory by an array, whose size will be the size of the set $A$. 

In a tree, each node has a parent node, except for the top node or root. In the array representation, each position of the array represents a node and the value of the node points to its parent node. Top nodes should have a value indicating that the node is the root of the tree. A negative number is used to signal the root, and the absolute value of this number is the number of elements of the block. 

A shallow tree is a tree in which all non-root nodes are connected directly to the root. To obtain a unique representation, we adopt the convention that the trees used in presenting a partition are shallow, and that the smallest value in a block will be the root node. A partition representation that respects the above conventions will be called a normalized (representation of a) partition. 

Thus, for $l,m \ge 1$, the array 
\[
\begin{array}{ccccccc}
[&\ldots&,-l,&\ldots        &,m,&\ldots &]\\
 &         &\uparrow&    &\uparrow& \\
 &         &\mbox{position $i$}& &\mbox{position $k$}  
\end{array}
\] means that $i$ is the root of a block (as its entry is negative) of size  $l$; in addition the element $k$ belongs to the block whose root is $m$.

Using these conventions, the encoding of the partition {\tiny\verb+[[1], [2], [3], [4], [5], [6]]+} induced by the identity relation is 

\begin{Verbatim}[fontsize=\tiny]
    [-1, -1, -1, -1, -1, -1]
\end{Verbatim}

\noindent while the partition {\tiny\verb+[[1, 2, 3, 4, 5, 6]]+} with just one block is encoded as
\begin{Verbatim}[fontsize=\tiny]
     [-6, 1, 1, 1, 1, 1]
\end{Verbatim}

Other examples are given below.
\begin{center}
\begin{tabular}{l|l}
\multicolumn{1}{c}{partition} & \multicolumn{1}{c}{encoded as} \\
\hline
{\tiny \verb+[[1, 6], [2], [3, 5], [4]]+} & {\tiny \verb+[-2, -1, -2, -1, 3, 1]+}\\
{\tiny \verb+[[1, 3, 5], [2, 6], [4]]+} & {\tiny \verb+[-3, -2, 1, -1, 1, 2]+}\\
{\tiny \verb+[[1, 2, 5, 6], [3, 4]]+} & {\tiny \verb+[-4, 1, -2, 3, 1, 1]+}\\
\end{tabular}
\end{center}
This non-intuitive way of representing partitions will prove its usefulness later on. 

\subsection{Congruences}
\label{subsection:congruences}

A congruence of an algebra is an equivalence relation on its underlying set that is compatible with all algebraic operations \cite{BURRIS}. 

Technically, a congruence relation is an equivalence relation $\equiv$ on an algebra that satisfies $\mu(a_1, a_2, ..., a_n)\equiv\mu(a'_1, a'_2, ..., a'_n)$ for every $n$-ary operation $\mu$ and all elements $a_1, ..., a_n, a'_1, ..., a'_n$ such that $a_i \equiv a'_i$ for each $i=1, ..., n$. 

Congruences will be represented in the same way as their corresponding partitions, as detailed in the previous section. 

\subsection{Homomorphism, Endomorphism, Isomorphism and Automorphism}
\label{subsection:morphisms}

 If $A$, $B$ are two algebras of the same type, then a function $f:A\to B$ is a \emph{homomorphism} from $A$ to $B$ if $\mu(f(a_1), f(a_2), ..., f(a_n))= f(\mu(a_1, a_2, ..., a_n))$ for every $n$-ary operation $\mu$ and $a_1, \dots a_n\in A$ (more precisely, $\mu$ here stands for the two operations of $A$ and $B$ that are indexed equally by the (common) index scheme). 

An \emph{endomorphism} is a \emph{homomorphism} from an algebra to itself, a \emph{monomorphism} is an injective \emph{homomorphism}, while an \emph{isomorphism} is a bijective \emph{homomorphism}. 
An \emph{automorphism} is a \emph{homomorphism} that is both an \emph{isomorphism} and an \emph{endomorphism}. 

\section{The Algorithms}
\label{algorithm}

\subsection{Congruences algorithm}
\label{subsection:congruencesAlg}

The starting point for the computation of congruences are the algorithms described in Freese \cite{FREESE} to calculate the smallest congruence containing a given partition $\Theta$ of a finite algebra $A$; in particular this is used to compute the smallest congruence containing a pair of elements $(a,b)\in A\times A$,  called the {\em principal congruence} generated by $\{a,b\}$. From this base algorithm all congruences of the algebra $A$ are generated in an efficient way.  
%
Optimizations to the original algorithm were introduced taking into account that we only allow operations of arity at most $2$. In addition, we take advantage of $C$ to get a faster implementation than the original implementations made by Freese and his collaborators. Finally, our implementation is integrated with GAP and hence fully compatible with its other resources. 

\subsubsection{Partition Functions}
\label{subsubsection:PartitionFunctions}

There are several partition functions that play an important role in the algorithm to compute principal congruences. 

The function \textbf{CreamRootBlock} is an operation that returns the root node for a node i. The root node of a node is itself if the value of the node representation is negative. If the value of the node representation  points to a different node then that node is the parent node. This algorithm could be run recursively until reaching the root node but is implemented iteratively for better performance. This operation will work both for normalized and non-normalized partitions. 

Before returning, the node parent of i is set to the found root node in order to make the tree representing the partition as shallow as possible, thus avoiding that in future calls the algorithm needs to transverse several nodes to reach the root node.

\begin{algorithm}[H]
\caption{CreamRootBlock (i, partition)}
\begin{algorithmic}
\STATE $j \leftarrow i$
\WHILE {$partition[j] \geq 0$}
    \STATE $j \leftarrow partition[j]$
\ENDWHILE
\IF{$i \neq j$}
    \STATE $partition[i] \leftarrow j$
\ENDIF
\RETURN $j$
\end{algorithmic}
\label{principalCongruence}
\end{algorithm}

The function \textbf{CreamJoinBlocks} is an operation that joins the blocks containing given elements x and y. This operation will work both for normalized and non-normalized partitions. The resulting partition may not in general be normalized  even if the original partition is.  

In order to keep the tree representing the partition as shallow as possible, the root node of the merged block will be the root node of the larger original block. 

\begin{algorithm}[H]
\caption{CreamJoinBlocks (x, y, partition)}
\begin{algorithmic}

\STATE $r \leftarrow CreamRootBlock (x,partition)$
\STATE $s \leftarrow CreamRootBlock (y,partition)$
\IF{$r \neq s$}
    \IF{$partition[r] < partition[s]$}
        \STATE $partition[r] \leftarrow partition[r] + partition[s]$
        \STATE $partition[s] \leftarrow r$
    \ELSE
        \STATE $partition[s] \leftarrow partition[r] + partition[s]$
        \STATE $partition[r] \leftarrow s$        
    \ENDIF
\ENDIF
\end{algorithmic}
\label{CreamJoinBlocks}
\end{algorithm}

The \textbf{CreamNumberOfBlocks} function returns the number of blocks of a partition. Given the encoding of partitions we simply count the number of positions of the array that have negative values.

\begin{algorithm}[H]
\caption{CreamNumberOfBlocks (partition)}
\begin{algorithmic}
\STATE $nblocks \leftarrow 0$
\STATE $dimension \leftarrow ArraySize(partition)$ 
\FOR{$i=1$ to $dimension$}
    \IF{$partition[i] < 0$}
        \STATE $nblocks \leftarrow nblocks + 1$
    \ENDIF    
\ENDFOR
\RETURN $nblocks$
\end{algorithmic}
\label{CreamNumberOfBlocks}
\end{algorithm}

The \textbf{CreamNormalizePartition} function normalizes the partition by making it shallow and having the smallest element of each block the root node.\\
\begin{algorithm}[H]
\caption{CreamNormalizePartition (partition)}
\begin{algorithmic}
\STATE $dimension \leftarrow ArraySize(partition)$
\FOR{$i=1$ to $dimension$}
    \STATE $r \leftarrow CreamRootBlock(i,partition)$
    \IF{$r \geq i$}
        \STATE $partition[i] \leftarrow -1$
        \IF{$r > i$}
            \STATE $partition[r] \leftarrow i$
        \ENDIF
    \ELSE
        \STATE $partition[r] \leftarrow partition[r] - 1$
    \ENDIF
\ENDFOR
\end{algorithmic}
\label{CreamNormalizePartition}
\end{algorithm}

The \textbf{CreamJoinPartition} function computes the join of two partitions, i.e. the smallest partition containing both input partitions. 


\begin{algorithm}[H]
\caption{CreamJoinPartition (partition1, partition2)}
\begin{algorithmic}
\STATE $dimension \leftarrow ArraySize(partition1)$
\FOR{$i=1$ to $dimension$}
        \STATE CreamJoinBlocks($i$,CreamRootBlock ($i$,$partition1$),$partition2$)
\ENDFOR
\end{algorithmic}
\label{JoinPartition}
\end{algorithm}
The \textbf{CreamComparePartitions} function compares normalized partitions. Returns 0 if the partitions are equal, -1 if partition1 $>$ partition2, and 1 if partition1 $<$ partition2. Partitions are ordered by the order of the underlying list. The purpose of this function is to allow binary search in sets of partitions.\\
\begin{algorithm}[H]
\caption{CreamComparePartitions (partition1, partition2)}
\begin{algorithmic}
\STATE $dimension \leftarrow Length(partition1)$
\FOR{$i=1$ to $dimension - 1$}
    \IF{$partition1[i] > partition2[i]$}
        \RETURN $-1$
    \ENDIF
    \IF{$partition1[i] < partition2[i]$}
        \RETURN $1$
    \ENDIF
\ENDFOR
\RETURN $0$
\end{algorithmic}
\label{ComparePartitions}
\end{algorithm}

\subsubsection{Computing Principal Congruences}
\label{subsubsection:principallCongAlg}

The base algorithm in \cite{FREESE} computes the smallest congruence containing $\Theta$, a partition of a finite algebra $A$. The simplest case is when $\Theta$ only contains one nontrivial block and this block has only two elements. 

The algorithm only works with unary operations. Since our algebras have binary operations, we have to convert them  to a family of  unary operations. In this case a binary operation $f(x,y)$ can be converted into $2n$ unary operations (where $n$ is the size of the algebra) by assigning to each element $x\in A$ the unary operations $c_x$ and $r_x$, which are induced by the corresponding column and row in the Cayley table of $f$:
$$c_x(y) = f(y,x) \mbox{ and }r_x(y) = f(x,y).$$
For example,  the algebra:
$$
[ [ [2, 1, 1], [1, 2, 2], 
    [1, 3, 2] ] ]
$$
induces the following  unary operations (after removing duplicates):
$$
[ [2, 1, 1], [1, 2, 2], [1, 3, 2], [1, 2, 3] ]. 
$$

The principal congruence algorithm takes as input the algebra and the generating pair of elements.

The algorithm joins the pair of elements in the same block and applies the algebra's unary functions to each element of the pair. In this way, it obtains one additional pair per function that will be joined in the same block, repeating this process until the list of pairs is exhausted.  

\begin{algorithm}[H]
\caption{CreamPrincipalCongruence (algebra,InitialPair)}
\begin{algorithmic}
\STATE $PairList \leftarrow [InitialPair]$
\STATE $partition \leftarrow SingletonPartition()$
\STATE $partition \leftarrow CreamJoinBlocks (partition, InitialPair[1], InitialPair[2])$
\STATE $NFuncs \leftarrow ArraySize(algebra)$
\WHILE {$PairList \neq empty$} 
    \STATE $Pair \leftarrow Pop (PairList)$
    \FOR{$i=1$ to $NFuncs$}
        \STATE $f \leftarrow algebra (i)$
        \STATE $r \leftarrow CreamRootBlock (partition, f(Pair[1]))$
        \STATE $s \leftarrow CreamRootBlock (partition, f(Pair[2]))$
        \IF{$r \neq s$}
            \STATE $partition \leftarrow CreamJoinBlocks (partition, r, s)$
            \STATE $PairList \leftarrow Push(PairList, [r,s])$
        \ENDIF    
    \ENDFOR
\ENDWHILE
\RETURN $partition$
\end{algorithmic}
\label{CreamPrincipalCongruence}
\end{algorithm}

As described in \cite{FREESE} this algorithm is very efficient showing a moderate growth of execution time with $n$. Apart from the algorithm itself, several implementation choices were used to increase efficiency as well.

One of these aspects is the use of arrays to encode partitions which allow for a constant and very fast random access to the partition element.

This is combined with the balancing and collapsing of the trees representing the blocks in the partition that are parts of  the \textbf{CreamJoinBlocks} and \textbf{CreamRootBlock} algorithms. This approach aims  for trees that are as shallow as possible during the execution of the \textbf{CreamPrincipalCongruence} algorithm.

In \textbf{CreamJoinBlocks} the tree is balanced by keeping as root of the joined block the root of the bigger original block.

In \textbf{CreamRootBlock} the tree is collapsed by setting the parent node of the element on which the function is called equal to its return value, avoiding having to traverse several nodes of the tree in future calls.

Both of these implementation details allow for shallower trees representing the blocks, which will make \textbf{CreamRootBlock} faster, given that fewer nodes will have to be transversed to determine the root node of a node.

This is especially important since \textbf{CreamRootBlock} is the most called function when calculating a principal congruence of an algebra.

To compute all the principal congruences this function is called for all pairs of elements of $A$.

\begin{algorithm}[H]
\caption{CreamAllPrincipalCongruences (algebra)}
\begin{algorithmic}
\STATE $dimension \leftarrow SizeAlgebra(algebra)$
\STATE $allPrincipalCongruences \leftarrow []$
\FOR{$i=1$ to $dimension-1$}
    \FOR{$j=i+1$ to $dimension$}
        \STATE $congruence \leftarrow CreamPrincipalCongruence (algebra, [i,j])$
        \STATE $AddCongruence (allPrincipalCongruences, congruence)$
    \ENDFOR
\ENDFOR
\RETURN $allPrincipalCongruences$
\end{algorithmic}
\label{CreamAllPrincipalCongruences}
\end{algorithm}

\subsubsection{Calculating All Congruences}
\label{subsubsection:AllCongAlg}

It is known that the congruences of an algebra form a lattice with the usual set inclusion partial order.  The meet is the usual intersection; the join of two congruences is the smallest congruence containing both. 

The minimal elements of this lattice are precisely the principal congruences, and we can obtain all congruences by computing joins, starting with the principal congruences.

The key part is to find an efficient way to combine minimal congruences to get all congruences of the algebra. The principal congruences are stored in an ordered set and the congruences are combined from start to end by the join operation. Each congruence resulting from these joins will be added to the ordered set (eliminating duplicates). 


\begin{algorithm}[H]
\caption{CreamAllCongruences (algebra)}
\begin{algorithmic}
\STATE $allPrincipalCongruences \leftarrow CreamAllPrincipalCongruences (algebra)$
\STATE $allCongruences \leftarrow allPrincipalCongruences$
\STATE $i \leftarrow 1$
\WHILE {$i < Size(allCongruences)$}
    \FOR{$j=1$ to $Size(allPrincipalCongruences)$}
        \IF{$\neg isContained (allPrincipalCongruences[j],allCongruences[i])$}
            \STATE $congruence \leftarrow JoinPartition (allCongruences[i], allPrincipalCongruences[j])$
            \STATE $AddCongruence (allCongruences, congruence)$
        \ENDIF    
    \ENDFOR
    \STATE $i \leftarrow i+1$
\ENDWHILE
\RETURN $allCongruences$
\end{algorithmic}
\label{CreamAllCongruences}
\end{algorithm}

The join operation is only called if the principal congruence is not contained in the congruence that is being joined with. 
This is optimized by preserving the information about which principal partitions were joined. 


\subsubsection{Calculating Monolithic Algebras}
\label{subsubsection:MonolithicAlg}

Monolithic Algebras are algebras that have a single minimal congruence that is contained in every other congruence except for the identity congruence.

This is simple to calculate having the list of minimal congruences and a function that compares two partitions and calculates whether one is contained in the other.

The \textbf{CreamContainedPartition} function returns whether a partition is contained in another.
This function is used internally to determine whether an Algebra is monolithic or not. This algorithm assumes that the partitions are normalized.

\begin{algorithm}[H]
\caption{CreamContainedPartition (partition1, partition2)}
\begin{algorithmic}
\STATE $dimension = ArraySize(partition1)$
\FOR{$i=1$ to $dimension$}
    \IF{$partition1[i] < 0$}
        \IF{$partition2[i] < 0$}
            \STATE $block2 \leftarrow i$
        \ELSE
            \STATE $block2 \leftarrow partition2[i]$
        \ENDIF
        \FOR{$j=i+1$ to $dimension$}
            \IF{$partition1[j] = i$ and not $block2 = partition2[j]$}
                \RETURN $false$
            \ENDIF
        \ENDFOR
    \ENDIF
\ENDFOR
\RETURN $true$
\end{algorithmic}
\label{CreamContainedPartition}
\end{algorithm}

We calculate the principal congruences and compare them eliminating those that contain other congruences. If a single congruence is returned then this congruence is contained in every other congruence (and the algebra is monolithic). This is done with the \textbf{CreamIsAlgebraMonolithic} function.

\begin{algorithm}[H]
\caption{CreamIsAlgebraMonolithic (algebra)}
\begin{algorithmic}
\STATE $partitions \leftarrow CreamAllPrincipalCongruences (algebra)$
\STATE $npartitions \leftarrow Size (partitions)$
\IF{$npartitions > 1$}
    \STATE $i \leftarrow 2$
    \REPEAT
        \IF{$CreamContainedPartition (partitions[1],partitions[i])$}
            \STATE $Remove (partitions,i)$
        \ELSIF{$CreamContainedPartition (partitions[i],partitions[1])$}
            \STATE $partitions[1] \leftarrow partitions[i]$
			\STATE $Remove (partitions,i)$
			\STATE $i \leftarrow 2$
		\ELSE
		    \STATE $i \leftarrow i + 1$
        \ENDIF		 
    \UNTIL{Size (partitions) $<$ i}
\ENDIF
\STATE $npartitions \leftarrow Size (partitions)$
\IF{$npartitions > 1$}
    \RETURN $false$
\ELSE
    \RETURN $true$
\ENDIF
\end{algorithmic}
\label{CreamIsAlgebraMonolithic}
\end{algorithm}

\subsection{Automorphism algorithm}
\label{subsection:automorphismAlg}
The automorphisms of an algebra $A$ of type ($2^m$,$1^n$) can be derived from the calculation of the automorphisms of a set of algebras each containing only one  binary operation of $A$. The automorphism group of the algebra will be the intersection of all the automorphism groups of these magmas that also commute, function composition-wise, with all the unary operations of the original algebra.  

%

In general, it is difficult to obtain automorphisms efficiently.  The authors of  the Loops GAP package \cite{LOOPS3_4_1} used an idea that  we could apply to general magmas. With some adaptations this allows for very effective computing of the automorphisms of magmas.
 
%
\subsubsection{Invariant Vector}
\label{subsubsection:invariant}
The general idea is to pick a list of properties invariant under homomorphisms and then partition the algebras using these properties. 
 For example, if $e$ is idempotent, then any endomorphism must map $e$ onto an idempotent.   This is the basis of the {\em Discriminator} in the Loops package, which implemented nine such invariants for the domain elements of quasigroups. However, most of these invariants cannot be carried over directly to magmas. We devised a total of  seventeen invariants that hold for any magma.   In general, unless otherwise said, given an element $x$ in a Magma $M$ and $k \in \mathbb{Z}^+$, we define $x^k$ to be $x^{k-1}*x$ for $k\ge 2$ and $x^1=x$ (ie. we associate on the left). For each element $p$ in a magma $M$, CREAM computes the following: 

\begin{enumerate}
	\item\label{1} Smallest $k$ such that $p^k = p^n$, with $n > k > 1$.
	\item\label{2} Number of domain elements for which $p$ is a left identity.
	\item\label{3} Number of domain elements for which $p$ is a right identity.
	\item\label{4} Number of elements $y$ such that $p = (py)p$. 
	\item\label{5} Number of distinct elements in row $p$ of the multiplication table.
	\item\label{6} Number of distinct elements in column $p$ of the multiplication table. 
	\item\label{7} 1 if $p$ is an idempotent, 0 otherwise.
	\item\label{8} Number of idempotents in column $p$ of the multiplication table.
	\item\label{9} Number of idempotents in row $p$ of the multiplication table.
	\item\label{10} 1 if the equality $p(pp) = (pp)p$ holds, 0 otherwise.
	\item\label{11} Number of domain elements that commute with $p$. 
	\item\label{12} Number of domain elements $s$ for which $(ss)p=p(ss)$.
	\item\label{13} Number of domain elements $s$ such that $s^2=p$.
	\item\label{14} Number of domain elements $s$ satisfying $p(ps) = (pp)s$
	\item\label{15} Number of multisets $\{x,y\}$ with $xy = yx = p$.
	\item\label{16} Number of elements $t\in M$ such that for two idempotents $e,f \in M$, we have $p=et=tf$.
	\item\label{17} Number of elements $t \in M$ for which there exist two elements $x,y \in M$ such that $p=xy$ and $t=yx$.
\end{enumerate}

Many of these invariants have an algebraic meaning. 
For example, (\ref{17}) is the size of the conjugacy class of $p$ when the magma is a group; (\ref{16}) is the size of the filter generated by $p$  when the magma is a regular semigroup;  (\ref{13}) is the number of square roots of $p$; (\ref{11}) is the size of the centralizer of  $p$; (\ref{8}) and (\ref{9}) are the number of idempotents in the left or right ideals generated by $p$; (\ref{5}) and (\ref{6}) are the number of elements in the left or right ideals generated by $p$; (\ref{1}) deals with the periodicity of $p$; etc. 

The intersection of the blocks induced by   these invariants partitions the elements of the algebra, and by a straightforward argument on first order logic any endomorphism must preserve the blocks of this partition. 
\subsubsection{Generating Set}
\label{subsubsection:generatingSet}
We can cut down the size of the search tree by focusing only on the images  of the elements in a generating set.
Some generating sets may be more suitable for finding automorphisms than others:

\begin{itemize}
	\item A smaller generating set is preferred because it would require fewer trials in constructing an isomorphism.
	\item Generators from partition blocks with fewer elements are preferred since any homomorphism must send the elements of a block into itself. 
\end{itemize}

\begin{samepage}
Our algorithm satisfies both the constraints above and efficiently produces a generating set. It is depicted in Algorithm \ref{efficentGSet}. 

\begin{algorithm}
\caption{Constructing Efficient Generating Set}
\begin{algorithmic}
	\STATE $submagma \leftarrow [\ ]$
	\STATE $generator \leftarrow [\ ]$
	\STATE $candidates \leftarrow $ magma elements sorted by ascending block size
	\WHILE{$submagma \neq magma$}
		\STATE $generator \leftarrow generator\ \cup$ the element that increase the size of $submagma$ most, and if 2 elements increase the size by the same amount, take the one from the smaller block.
		\STATE $submagma \leftarrow $ submagma generated by $generator$
		\STATE $candidates \leftarrow candidates$ with elements in submagma removed
	\ENDWHILE
\end{algorithmic}
\label{efficentGSet}
\end{algorithm}
\end{samepage}
\subsubsection{The Automorphism Algorithm}
\label{subsubsection:automorphismAlgS}

To find the automorphisms we start by partitioning the Magma according to the invariant vectors of each element and then obtain an efficient generating set. 
We then need to find the images of  the generating set taking into account that  each element in the generating set can only map to elements having the same invariant vector. Then we extend this  partial map to a full map and check that the homomorphism condition is satisfied.

We tested the algorithm on groups, loops, and quasigroups, obtaining the same results as yielded by the existing tools. For magmas of orders 2 and 3, and semigroups of orders 6 and 7 (as there are no general tools to compute automorphisms) we double check  our results  against the output of Mace4 \cite{MACE}. 


\subsubsection{Algebras of Type ($2^m$,$1^n$)}
\label{subsubsection:algebra_type_2}

Given an algebra of type ($2^m$,$1^n$) its automorphism group can be computed by taking the intersection of all the automorphism groups of the binary operations that also commute, function composition-wise, with all the unary operations. 

%

A useful observation is that the intersection of the trivial group and any group is the trivial group. Therefore, if any of the automorphism group generated in the process is a trivial group, then the search can be terminated with the trivial group declared the automorphism group for the algebra. 

This algorithm is implemented in the \textbf{CreamAutomorphisms} function in the CREAM package, as shown in Algorithm~\ref{CreamAutomorphisms}.

\begin{algorithm}[H]
	\caption{CreamAutomorphisms(algebra)}
	\begin{algorithmic}
		\STATE $binaryOperations \leftarrow AllBinaryOperations(algebra)$
		\STATE $unaryOperations \leftarrow AllUnaryOperations(algebra)$
		\STATE $D \leftarrow \emptyset$
		\REPEAT
			\STATE $B \leftarrow$ Pop a binary operation from $binaryOperations$
			\STATE $C$ $\leftarrow$ all automorphisms of $B$
			\STATE $C$ $\leftarrow$ all automorphisms in $C$ that commute with every unary operation in $unaryOperations$
			\IF {$C=\emptyset$}
				\RETURN $\emptyset$
			\ENDIF
			\STATE $D \leftarrow{D}\cup{C}$
		\UNTIL{$binaryOperations=\emptyset$}
		\RETURN Intersection($D$)
	\end{algorithmic}
	\label{CreamAutomorphisms}
\end{algorithm}

\subsection{Endomorphisms algorithm}
\label{subsection:endomorphismsAlg}
The classic approach to calculate endomorphisms would be to use \textbf{Mace4} to search for endomorphisms and while this is very effective and fast for low order algebra further optimizations are needed for high order algebras in which the congruences of the algebra can be used to limit the \textbf{Mace4} search space for endomorphisms.

To compute all endomorphisms of high order algebras, CREAM does the following: 
\begin{enumerate}
    \item Compute the congruences of the algebra $A$;
    \item For each congruence $R$, compute $A/R$ (except for the trivial congruence);
    \item Compute the subalgebras of $A$ that are isomorphic to $A/R$ (using the finite model builder Mace4);
    \item For each compatible pair subalgebra/congruence derive a corresponding endomorphism;
    \item Add the automorphisms of $A$.
\end{enumerate}

The congruences of $A$ are calculated with the algorithms described in \ref{subsection:congruencesAlg}. 
For each congruence $R$ (except for the trivial congruence that will be dealt  with later), a representation of the $A/R$ operation can be obtained by replacing the values in $A$'s operation tables by their roots, and eliminating rows and columns not corresponding to roots. If we have the algebra:\\

\begin{Verbatim}[fontsize=\tiny]
[ [ [ 1, 1, 1, 1, 1, 1 ], [ 1, 1, 1, 1, 1, 1 ], [ 1, 1, 1, 1, 1, 1 ], 
    [ 1, 1, 1, 1, 1, 1 ], [ 1, 1, 1, 1, 1, 1 ],[ 1, 1, 1, 2, 3, 1 ] ] ] 
\end{Verbatim}
and the congruence
\begin{Verbatim}[fontsize=\tiny]
[ [ 1, 2, 3], [4], [5, 6] ] - [ -3, 1, 1, -1, -2, 5 ],
\end{Verbatim}
then the  $A/R$ operation is: 
{\tiny
\[
\begin{tabular}{c|ccc}
  &1&4&5\\ \hline
1&1&1&1\\ 
4&1&1&1\\ 
5&1&1&1\\ 
\end{tabular}
\]
%
}
and passed to Mace4 as follows:
\begin{Verbatim}[fontsize=\tiny]
f(a0,a0)=a0.
a0!=a3.
f(a0,a3)=a0.
a0!=a4.
f(a0,a4)=a0.
f(a3,a0)=a0.
f(a3,a3)=a0.
a3!=a4.
f(a3,a4)=a0.
f(a4,a0)=a0.
f(a4,a3)=a0.
f(a4,a4)=a0.
\end{Verbatim}
where $f$ is the binary operation of the algebra. It is worth observing  that Mace4 is zero-based, unlike GAP, which is one-based.

Running Mace4 with the definition of the algebra operation and the above encoding of $A/R$ as assumptions we get 16 possible models each one of them corresponding to a different endomorphism of the algebra $A$. One of these models would be in Mace4's output:

\begin{Verbatim}[fontsize=\tiny]
interpretation( 6, [number = 1,seconds = 0], [
    function(a0, [0]),
    function(a3, [1]),
    function(a4, [2]),
    function(fa1(_,_), [
        0,0,0,0,0,0,
        0,0,0,0,0,0,
        0,0,0,0,0,0,
        0,0,0,0,0,0,
        0,0,0,0,0,0,
        0,0,0,1,2,0])]).
\end{Verbatim}

From this model we get  the partial mapping $f:\{1,4,5\}\subseteq A\to \{1,2,3\}$, defined by $f(1)=1$, $f(4)=2$ and $f(5)=3$, that can be  represented by $[ 1, , , 2, 3, ]$. 
The gaps can be easily filled since elements of the same block must map to the same image. In this example, $2$ and $3$ have the same image as $1$, while $6$ and  $5$ have the same image yielding the endomorphism: $
[ 1, 1, 1, 2, 3, 3]$. 

Repeating this process for all congruences and for all models obtained with Mace4 from each congruence, all the algebra's endomorphisms (except for the automorphisms) can be obtained.

Finally, for the trivial congruence the associated endomorphisms would be also automorphisms and in this case it is more efficient to use the the algorithm described in \ref{subsection:automorphismAlg} than using Mace4.

\subsection{Monomorphisms algorithms}
\label{subsection:monomorphismsAlg}

Invariants can be used to speed up the process of finding an monomorphism from one magma to another, or from one algebra to another.  If $f$ is an injective homomorphism from a magma $A$ to a magma $B$, then it is a isomorphism from $A$ to $B$ restricted to the range of $f$.  Thus, the same ideas of applying invariants in constructing automorphisms can be used for constructing monomorphisms.  Specifically, an monomorphism $f$ can map an element $a\in{A}$ to $b\in{B}$ only if $a$ and $b$ have the same invariant vector.  This greatly reduces the search space for monomorphisms between $A$ and $B$.

\begin{algorithm}[H]
	\caption{CreamAllMonomorphismMagmas(magma1, magma2)}
	\begin{algorithmic}
		\STATE $genL \gets$ generating set from $magma1$
		\FOR{$i$ = 1 to $Size(magma2)$}
			\STATE $rangeM[i] \gets$ list of elements of $magma2$ that have same invariant vector as $genL[i]$
		\ENDFOR
		\STATE $monoMaps \gets \emptyset$
		\FORALL{mappings $m$ from $genL$ to elements of $rangeM$ s.t. $genL[i]$ maps only to elements in $rangeM[i]$}
			\IF{$m$ is an injective homomorphism} 
				\STATE $monoMaps \gets monoMaps \cup m$
			\ENDIF
		\ENDFOR
		\RETURN $monoMaps$
	\end{algorithmic}
	\label{CreamAllMonomorphismMagmas}
\end{algorithm}

\begin{algorithm}[H]
	\caption{CreamAllMonomorphism(algebra1, algebra2)}
	\begin{algorithmic}
		\STATE b1 $\leftarrow$ first binary operation in $algebra1$
		\STATE b2 $\leftarrow$ first binary operation in $algebra2$
		\STATE $monoList$ $\leftarrow$ CreamAllMonomorphismMagmas($b1, b2$)
		\STATE $algebraMonoList \leftarrow []$ 
		\FOR{each $monoMap$ in $monoList$}
			\STATE $mono \leftarrow true$
			\FOR{each corresponding pair of operations $b1\in{algebra1}, b2\in{algebra2}$}
				\IF{$monoMap$ is not an isomorphism from $b1$ to $b2$}
					\STATE $mono \leftarrow false$ \\
					break out of for loop
				\ENDIF
			\ENDFOR
			\IF{$mono$}
			    \STATE Add $monoMap$ to $algebraMonoList$ if not already in the list
			\ENDIF
		\ENDFOR
		\RETURN $algebraMonoList$
	\end{algorithmic}
	\label{CreamAllMonomorphisms}
\end{algorithm}

In this section only the function CreamAllMonomorphisms is addressed but the Cream package also includes the functions CreamExistsMonomorphism (returning true if a monomorphism exists) and CreamOneMonomorphism (returning one monomorphism).

\subsection{Epimorphisms algorithm}
\label{subsection:epimorphismsAlg}

CreamAllEpimorphisms($A_1$, $A_2$) returns all epimorphisms from $A_1$ to $A_2$, provided the algebras are compatible. The algorithm first finds all congruences $\varphi$ of $A_1$ such that $A_1/\varphi$  and $A_2$ are isomorphic. For each such $\varphi$, the corresponding epimorphisms are obtained by composing the quotient map with an isomorphism to $A_2$ and all automorphisms of $A_2$. 

For efficiency, a size check is implemented before searching for isomorphisms  and the automorphisms of $A_2$ are only calculated once.


\begin{algorithm}[H]
	\caption{CreamAllEpimorphisms(algebra1, algebra2)}
	\begin{algorithmic}
	    \STATE $allCongruences \leftarrow CreamAllCongruences(algebra1)$
	    \STATE $dimension1 \leftarrow SizeAlgebra(algebra1)$
	    \STATE $dimension2 \leftarrow SizeAlgebra(algebra2)$
	    \STATE $autoList \leftarrow []$
	    \STATE $epiList \leftarrow []$
	    \FORALL {$cong$ in $allConguences$}
            \IF {$dimension2 = CreamNumberOfBlocks (cong)$}       
                \STATE $[qalgebra,mapToQalgebra] \leftarrow QuotientAlgebraFromCongruence (algebra1, cong)$
                \STATE $iso \leftarrow IsomorphismAlgebras (qalgebra,algebra2)$
                \IF {not $iso = fail$}
                    \IF {autoList = []}
                        \STATE $autoList \leftarrow CreamAutomorphisms (algebra2)$
                    \ENDIF
                    \FORALL {$auto$ in $autoList$}
                        \STATE $epi \leftarrow []$
                        \FOR {$j=1$ to $dimension1$}
                            \STATE $epi[j] \leftarrow auto[iso[mapToQalgebra[i]]]$
                        \ENDFOR
                        \STATE Add $epi$ to $epiList$ if not already in the list
                    \ENDFOR
                \ENDIF
            \ENDIF
	    \ENDFOR
		\RETURN $epiList$
	\end{algorithmic}
	\label{CreamAllEpimorphisms}
\end{algorithm}

In this section only the function CreamAllEpimorphisms is addressed but the Cream package also includes the functions CreamExistsEpimorphism (returning true if a epimorphism exists) and CreamOneEpimorphism (returning one epimorphism).

\subsection{SubUniverses algorithm}
\label{subsection:subuniversesAlg}

The SubUniverses algorithm returns a list of all underlying sets of all subalgebras of the input algebra. It first generates all 1-generated subalgebras, then iteratively expands each $i$-generated algebra by adding another generator.\\

For performance enhancement, these expansions are limited to one  element from each orbit of the algebra's automorphism group. At the end of each cycle, the remaining $(i+1)$-generated subuniverses are obtained by applying the automorphism group. The algorithm stops when an iteration produces only one subuniverse with the same size of the algebra.\\

The algorithm relies on the  SubUniverseFromElement routine, which calculates the subalgebra generated by a subalgebra and an additional element. This routine assumes that the input subuniverse is closed under the algebra's operations. The algorithm adds the element to the updated subuniverse, calculates all directly generated additional elements, and places them in a temporary list. These elements are then filtered for those that are already in the subuniverse, and the remaining elements are put into a second list. Elements from the second list are then added to the subuniverse iteratively.  \\


\begin{algorithm}[H]
	\caption{CreamAllSubUniverses(algebra)}
	\begin{algorithmic}
	    \STATE $dimension \leftarrow CreamSizeAlgebra(algebra)$
	    \STATE $autoList \leftarrow CreamAutomorphisms (algebra)$
	    \FOR{$i=1$ to $dimension$}
	        \STATE $sigma[i] \leftarrow []$
	        \IF {$i=1$}
	            \STATE $sigmaMinus \leftarrow [[]]$
	        \ELSE 
	            \STATE $sigmaMinus \leftarrow sigmaExpanded[i-1]$
	        \ENDIF
	        \FORALL{$cSigma$ in $sigmaMinus$}
	            \STATE $elemList \leftarrow [1 .. dimension]$
	            \STATE Remove all elements in $cSigma$ from $elemList$
        	    \WHILE {$Size(elemList) <> 0$}
        	        \STATE $j \leftarrow elemList[1]$
        	        \STATE $sUniverse \leftarrow SubUniverseFromElement (algebra, cSigma, j)$
        	        \STATE Add $sUniverse$ to $sigma[i]$ if not already in the list
        		    \STATE $autoEs \leftarrow $ the orbit of $j$ under $autoList$
        		    \STATE Remove all elements in $autoEs$ from $elemList$
        	    \ENDWHILE	            
	        \ENDFOR
	        \STATE $sigmaExpanded[i] \leftarrow []$
    	    \FORALL{$sUniverse$ in $sigma[i]$}
    	        \STATE $autoUniverses \leftarrow$ the orbit of $sUniverse$ under $autoList$
    	        \STATE Add all elements of $autoUniverses$ to $sigmaExpanded[i]$ if not already in the list
    	    \ENDFOR
    	    \IF {$Size(sigma[i]) = 1$ and $Size(sigma[i][1]) = dimension$}
    	        \STATE $break$
    	    \ENDIF
	    \ENDFOR	    
	    \STATE $sUniverses \leftarrow \cup_{j=1}^i sigmaExpended[j]$
    	\RETURN $sUniverses$
	\end{algorithmic}
	\label{CreamAllSubUniverses}
\end{algorithm}

\begin{algorithm}[H]
	\caption{SubUniverseFromElement(algebra, initSubUniverse, element)}
	\begin{algorithmic}
	    \STATE $dimension \leftarrow CreamSizeAlgebra(algebra)$
        \STATE $elemList \leftarrow [element]$	
        \STATE $newSubUniverse \leftarrow initSubUniverse$
        \STATE $nElem \leftarrow Size(elemList)$
        \WHILE {$nElem > 0$}
            \STATE $currentElem \leftarrow elemList[nElem]$
            \STATE Remove element in position $nElem$ from  $elemList$  
            \STATE $newElemList \leftarrow []$
            \FORALL {operations $op$ of the $algebra$}
                \IF {$op$ is binary}
                    \STATE Add $op[currentElem][currentElem]$ to $newElemList$ if not already in the list
                    \FORALL {$subUniverseElem$ in $newSubUniverse$}
                        \STATE Add $op[currentElem][subUniverseElem]$ to $newElemList$ if not already in the list
                        \STATE Add $op[subUniverseElem][currentElem]$ to $newElemList$ if not already in the list
                    \ENDFOR
                \ELSE
                    \STATE Add $op[currentElem]$ to $newElemList$ if not already in the list
                \ENDIF
            \ENDFOR
            \STATE Add $currentElem$ to $newSubUniverse$ if not already in the list   
            \STATE $elemList \leftarrow elemList \cup (newElemList \setminus newSubUniverse$)         
            \STATE $nElem \leftarrow Size(elements)$
            \IF {$Size(newSubUniverse) + nElem = dimension$}
                \STATE $newSubUniverse \leftarrow [1 .. dimension]$
                \STATE $nElem \leftarrow 0$
            \ENDIF
        \ENDWHILE
        \RETURN $newSubUniverse$
	\end{algorithmic}
	\label{SubUniverseFromElement}
\end{algorithm}

%
%
%
%
%

\subsection{DivisorUniverses algorithm}
\label{subsection:divisorsAlg}

CreamAllDivisorUniverses($A_1$, $A_2$) checks if $A_1$ has a divisor that is isomorphic to $A_2$, where $A_1$ and $A_2$ are compatible algebras. Its return is a list of all pairs $(B,\varphi)$, such that $B$ is a subuniverse of $A_1$, $\varphi$ is a congruence of the subalgebra $B$, and $B/\varphi$ is isomorphic to $A_2$.\\

The algorithm first calculates the subalgebras of $A_1$, then their quotients and then checks those for isomorphisms to $A_2$, in each case using the corresponding algorithms of the Cream package. A size condition is used to prune unnecessary calculations. 

\begin{algorithm}[H]
	\caption{CreamAllDivisorUniverses(algebra1, algebra2)}
	\begin{algorithmic}
	    \STATE $dimension \leftarrow CreamSizeAlgebra(algebra2)$
	    \STATE $subUniverseList \leftarrow CreamSubUniverses (algebra1)$
	    \STATE $divisorList \leftarrow []$
	    \FORALL{$subUniverse$ in $subUniverseList$}
	        \IF {$dimension <= Size (subUniverse)$}
	            \STATE $subAlgebra \leftarrow CreamSubUniverse2Algebra (algebra1, subUniverse)$
	            \STATE $congList \leftarrow CreamAllCongruences (subAlgebra)$
	            \STATE $rcongList \leftarrow []$
	            \FORALL{$cong$ in $congList$}
	                \IF{$dimension = CreamNumberOfBlocks (cong)$}
	                    \STATE Add $cong$ to $rcongList$ if not already in the list
	                \ENDIF
	            \ENDFOR
	            \FORALL{$cong$ in $rcongList$}
	                \STATE $[qAlgebra,mapToQalgebra] \leftarrow QuotientAlgebraFromCongruence (subAlgebra, cong)$
	                \IF{$CreamAreAlgebrasIsomorphic (qAlgebra, algebra2)$}
	                    \STATE Add $[subUniverse, cong]$ to $divisorList$ if not already in the list
	                \ENDIF
	            \ENDFOR
	        \ENDIF
    	\ENDFOR
    	\RETURN $divisorList$
	\end{algorithmic}
	\label{CreamAllDivisorUniverses}
\end{algorithm}

In this section only the function CreamAllDivisorUniverses is addressed but the Cream package also includes the functions CreamExistsDivisor (returning true if a divisor between the algebras exists) and CreamOneDivisorUniverse (returning one divisor between the algebras).

\subsection{DirectlyReducible algorithm}
\label{subsection:directlyReducibleAlg}

CreamAllDirectlyReducible($A$) looks for two (non-trivial)  commuting congruences $\phi, \psi$ of $A$, whose meet is trivial and whose join is $A \times A$. If those exist, $A$ is isomorphic to the direct product $A/\phi \times A/\psi$. \\

Because all algebras are finite, it suffices to check that the congruences have a trivial meet, and that the corresponding product algebra has the correct size. The algorithm uses numerical constraints to avoid unnecessary calculations.\\

The return of the algorithm is a list of pairs of congruences. These pairs should be consider as sets, i.e. if $(\alpha, \beta)$ is listed the (equally valid) pair $(\beta,\alpha)$ is not listed separately.

\begin{algorithm}[H]
	\caption{CreamAllDirectlyReducible(algebra)}
	\begin{algorithmic}
	    \STATE $dimension \leftarrow CreamSizeAlgebra(algebra)$
	    \STATE $congList \leftarrow CreamAllCongruences (algebra)$
	    \STATE $pCongList \leftarrow []$
	    \FORALL{$cong$ in $congList$}
	        \STATE $cNBlocks \leftarrow CreamNumberOfBlocks (cong)$
	        \IF {$cNBlocks <> 1$ and $cNBlocks <> dimension$ and $IsInt(dimension/cNBlocks)$}
	            \STATE Add $cong$ to $pCongList$
	        \ENDIF
    	\ENDFOR
    	\STATE $pairs \leftarrow []$
    	\FORALL {$2$-element subsets $\{i,j\}$ of $pCongsList$}
        	    \STATE $iNBlocks \leftarrow CreamNumberOfBlocks(i)$
        	    \STATE $jNBlocks \leftarrow CreamNumberOfBlocks(j)$
        	    \IF {$iNBlocks*jNBlocks=dimension$}
        	        \STATE $cMeet \leftarrow MeetPartition (i, j)$
        	        \STATE $cMeetNBlocks \leftarrow CreamNumberOfBlocks (cMeet)$
        	        \IF {$cMeetNBlocks = dimension$}
        	            \STATE Add $[i,j]$ to $pairs$
        	        \ENDIF
        	    \ENDIF
    	\ENDFOR
    	\RETURN $pairs$
	\end{algorithmic}
	\label{CreamAllDirectlyReducible}
\end{algorithm}

In this section only the function CreamAllDirectlyReducible is addressed but the Cream package also includes the functions CreamExistsDirectlyReducible (returning true if there are commuting congruences R, T where A/R x A/T is isomorphic to A) and CreamOneDirectlyReducible (returning one pair of commuting congruences R, T where A/R x A/T is isomorphic to A).

\section{Performance}
\label{performance}

\subsection{Congruences algorithm performance}
\label{subsection:congruencesAlgPerf}

The \textbf{semigroup} GAP package \cite{semigroups} includes the function \textbf{CongruencesOfSemigroup} that returns the congruences of a semigroup, for some classes of semigroups. For Rees matrix semigroups the function  \textbf{CongruencesOfSemigroup} takes advantage of theoretical results and hence achieves better results than \textbf{CreamAllCongruences}. But to achieve this kind of performance the semigroup needs to be created using the GAP functions \textbf{ReesMatrixSemigroup} or \textbf{ReesZeroMatrixSemigroup}. If instead  an isomorphic copy of a Rees Matrix (Zero) semigroup is created using a multiplication table  the performance will be much worse.

\begin{table}
	\caption{Performance Comparison on calculating all congruences of Rees Matrix  Semigroups.}\label{reesCongruencePerformance}
	\renewcommand\arraystretch{1.5}
	\noindent\[
{\tiny
	\begin{tabular}{|c|c|c|c|c|}
	\hline
    Size&CongruencesOfSemigroup&CongruencesOfSemigroup&CreamAllCongruences&UACalc\\
    &(ReesMatrixSemigroup)&(MultiplicationTable)&&\\
    &(ms)&(ms)&(ms)&(ms)\\
    \hhline{|=|=|=|=|=|}
    12&1&15&$<1$&1\\
    \hline
    18&2&28&$<1$&3\\
    \hline
    24&5&55&2&7\\
    \hline
    30&1&92&4&32\\
    \hline
    36&4&139&10&34\\
    \hline
    42&11&203&16&55\\
    \hline
    48&8&296&24&77\\
    \hline
    54&5&416&46&138\\
    \hline
    60&5&553&60&220\\
    \hline
    66&1&671&82&291\\
    \hline
    72&9&873&136&409\\
    \hline
    78&2&1 098&180&555\\
    \hline
    84&5&1 341&233&743\\
    \hline
    90&4&1 622&307&961\\
	\hline
	\end{tabular}
}
	\]
\end{table}

For Rees Matrix Semigroups or Rees Zero Matrix Semigroups (generated with \textbf{ReesMatrixSemigroup} or \textbf{ReesZeroMatrixSemigroup}) the function \textbf{CongruencesOfSemigroup} takes advantage of their particular structure by applying the efficient linked triple algorithm \cite{How95}. However, this efficiency can only be realized in this restricted setting. 

\textbf{UACalc} also implements Freese's algorithm, and was done under his supervision in Java and Jython. It includes a GUI version in Java and a command line interface in Jython. The performance of the GUI version is poor and it will not be used for comparison. On the other hand the command line interface has a performance that is comparable to \textbf{CreamAllCongruences}. The same random Rees Matrix Semigroups that were used with \textbf{CongruencesOfSemigroup} and \textbf{CreamAllCongruences} were written into files readable by \textbf{UACalc} and the congruences were calculated in \textbf{UACalc}.

\textbf{CreamAllCongruences} is consistently more than 3 times faster than \textbf{UACalc}. For these algebras with one binary operation the runtime rises roughly proportionally to $n^4$ (where $n$ is the size of the algebra) or $t^2$ (being $t$ the number cells of the multiplication matrix of the algebra operation). 

\subsection{Automorphisms algorithm performance}
\label{subsection:automorphismsAlgPerf}

The CREAM package automorphism function was run on algebras of Type ($2^m$,$1^n$) to compare the timings against the Loops package.  
All experiments were run 7 times, with the lowest and highest times discarded. The averages of the remaining 5 are reported in Table \ref{magmautvsloops}.  
We see that the speeds from both packages are quite close. 

\begin{table}[H]
	\caption{Performance Comparison between Loops and CREAM on Automorphism Group Generation.}\label{magmautvsloops}
	\renewcommand\arraystretch{1.5}
	\noindent\[
{\tiny	\begin{tabular}{|c|c|c|}
	\hline
	Algebraic Structure&Loops Package(sec)&CREAM Package(sec)\\
	\hhline{|=|=|=|}
	Quasigroups, order 5 & 0.588 & 0.562\\
	\hline
	Quasigroups, order 6 & 520 & 541\\
	\hhline{|=|=|=|}
	Loops, order 6 & 0.101 & 0.086\\
	\hline
	Loops, order 7 & 14.623 & 13.945\\
	\hhline{|=|=|=|}
	Groups, order 32 & 0.854 & 0.935\\
	\hline
	Groups, order 64 & 39.98 & 41.05\\
	\hline
	Groups, order 128 & 12,031 & 12,160\\
	\hline
	\end{tabular}
}	\]
\end{table}

\subsection{Endomorphisms algorithm performance}
\label{subsection:endomorphismsAlgPerf}

To evaluate the performance of the algorithms to calculate endomorphisms,  we  tested both the classic algorithm and the congruence algorithm running with \textbf{Mace4}, as described in \ref{subsection:endomorphismsAlg}, on several types of Semigroups, Monoids and Groups. The only GAP function to calculate endomorphisms that came to our attention is the function Endomorphisms from the package SONATA. This function can calculate endomorphisms for a very narrow set of algebras, namely groups and near-rings. The groups in the list algebras were also run with SONATA.\\ 

\begin{table}[H]
	\caption{Performance Comparison between endomorphisms calculation algorithms}\label{endo performance}
	\renewcommand\arraystretch{1.5}
	\noindent\[
{\tiny
	\begin{tabular}{|c|c|c|c|c|c|}
	\hline
    Algebra&Size&Number of&Classic&Congruences&SONATA\\
    &&Endomorphisms&(ms)&(ms)&(ms)\\
    \hhline{|=|=|=|=|=|=|}
    GossipMonoid(3)&11&66&30&1 635&NA\\
    \hline
    PlanarPartitionMonoid(2)&14&72&56&233&NA\\
    \hline
    JonesMonoid(4)&14&72&45&232&NA\\
    \hline
    BrauerMonoid(3)&15&28&38&160&NA\\
    \hline
    PartitionMonoid(2)&15&89&50&313&NA\\
    \hline
    FullPBRMonoid(1)&16&1 426&585&5 134&NA\\
    \hline
    SymmetricGroup(4)&24&58&101&142&166\\
    \hline
    FullTransformationSemigroup(3)&27&40&116&364&NA\\
    \hline
    FullTransformationMonoid(3)&27&40&112&360&NA\\
    \hline
    SymmetricInverseMonoid(3)&34&54&274&504&NA\\
    \hline
    JonesMonoid(5)&42&113&746&777&NA\\
    \hline
    MotzkinMonoid(3)&51&98&1 972&2 400&NA\\
    \hline
    PartialTransformationMonoid(3)&64&138&2 586&1 963&NA\\
    \hline
    PartialBrauerMonoid(3)&76&165&7 897&7 796&NA\\
    \hline
    BrauerMonoid(4)&105&274&50 164&19 875&NA\\
    \hline
    SymmetricGroup(5)&120&146&32 861&2 979&201\\
    \hline
    PlanarPartitionMonoid(3)&132&393&95 460&25 803&NA\\
    \hline
    JonesMonoid(6)&132&393&131 472&22 889&NA\\
    \hline
    PartitionMonoid(3)&203&687&741 441&105 421&NA\\
    \hline
    SymmetricInverseMonoid(4)&209&282&470 901&75 940&NA\\
	\hline
	\end{tabular}
}
	\]
\end{table}

In this test the classic algorithm not using congruences is faster than using congruences for algebras up to size 60 but for larger algebras its runtime raises very fast and the algorithm using congruences becomes faster.\\

In view of these results, the \textbf{CreamEndomorphisms} function uses the classic algorithm for algebras with size up to 60 and the congruences algorithm for larger algebras.\\

The SONATA package achieves very good results for groups but it only works on a very narrow set of algebras.\\

\section{Conclusion}
\label{conclusion}

The CREAM package provides efficient algorithms for the calculation of congruences, automorphisms and endomorphisms. These algorithms work with algebras of type $(2^m, 1^n)$, covering a very wide set of algebras and being consistently fast.


The integration with Mace4 allowed us to combine efficient algorithms such as \cite{FREESE} with the wide set of possibilities provided by a first order axiom-based model searcher. 

The appendices of this article include a detailed performance discussion of the main algorithms and applications of the CREAM package to different algebra classes.

The CREAM package can be found and downloaded from \url{https://gitlab.com/rmbper/cream} and the outputs from the test runs described in section \ref{section:applications} can be found in \url{https://gitlab.com/rmbper/cream_data}.

\clearpage

\begin{appendices}

\section{Detailed Performance}
\label{section:dperformance}

\subsection{Congruences algorithm performance}
\label{subsection:congruencesAlgPerfAnnex}

The \textbf{semigroup} GAP package \cite{semigroups} includes the function \textbf{CongruencesOfSemigroup} that returns the congruences of a semigroup. This function doesn't work for all semigroups, but for those in which it works  
we have run the \textbf{CongruencesOfSemigroup} function against our \textbf{CreamAllCongruences} 5 times and calculated the average runtime after removing the best and worst runtimes getting the following results:
\begin{table}[h]
	\caption{Performance Comparison between \textbf{CongruencesOfSemigroup} and CREAM on calculating all congruences.}\label{creamvssemigroups}
	\renewcommand\arraystretch{1.5}
	\noindent\[
	{\tiny
	\begin{tabular}{|c|c|c|c|}
	\hline
	Size&Number of Algebras&CongruencesOfSemigroup(ms)&CreamAllCongruences(ms)\\
	\hhline{|=|=|=|=|}
	1&1&7&$<1$\\
	\hline
	2&3&20&$<1$\\
	\hline
	3&4&29&$<1$\\
	\hline
	4&7&41&$<1$\\
	\hline
	5&9&53&$<1$\\
	\hline
	6&8&71&$<1$\\
	\hline
	7&12&153&4\\
	\hline
	8&14&484&35\\
	\hline
	1-8&58&858&40\\
	\hline
	\end{tabular}
	}
	\]
\end{table}

The number of semigroups for which \textbf{CongruencesOfSemigroup} works is very limited, yet even for those algebras \textbf{CreamAllCongruences} is more than 20 times faster. 

The  tests show that there is a very specific type of semigroup in which the function  \textbf{CongruencesOfSemigroup} takes advantage of theoretical results and hence achieves better results than \textbf{CreamAllCongruences}. These are the Rees Matrix Semigroups or the Rees Zero Matrix Semigroups. But to achieve this kind of performance the semigroup needs to be created using the GAP functions \textbf{ReesMatrixSemigroup} or \textbf{ReesZeroMatrixSemigroup}. If instead  an isomorphic copy of a Rees Matrix (Zero) semigroup is created using a multiplication table  the performance will be similar to what was found above. See in the following table the average runtime for 10 runs with random Rees Matrix Semigroups:
\begin{table}[h]
	\caption{Performance Comparison on calculating all congruences of Rees Matrix  Semigroup.}\label{reesAllCongPerformance}
	\renewcommand\arraystretch{1.5}
	\noindent\[
{\tiny
	\begin{tabular}{|c|c|c|c|c|}
	\hline
    Size&CongruencesOfSemigroup&CongruencesOfSemigroup&CreamAllCongruences&UACalc\\
    &(ReesMatrixSemigroup)&(MultiplicationTable)&&\\
    &(ms)&(ms)&(ms)&(ms)\\
    \hhline{|=|=|=|=|=|}
    12&1&15&$<1$&1\\
    \hline
    18&2&28&$<1$&3\\
    \hline
    24&5&55&2&7\\
    \hline
    30&1&92&4&32\\
    \hline
    36&4&139&10&34\\
    \hline
    42&11&203&16&55\\
    \hline
    48&8&296&24&77\\
    \hline
    54&5&416&46&138\\
    \hline
    60&5&553&60&220\\
    \hline
    66&1&671&82&291\\
    \hline
    72&9&873&136&409\\
    \hline
    78&2&1 098&180&555\\
    \hline
    84&5&1 341&233&743\\
    \hline
    90&4&1 622&307&961\\
	\hline
	\end{tabular}
}
	\]
\end{table}
For Rees Matrix Semigroups or Rees Zero Matrix Semigroups (generated with \textbf{ReesMatrixSemigroup} or \textbf{ReesZeroMatrixSemigroup}) the function \textbf{CongruencesOfSemigroup} takes advantage of their particular structure by applying the efficient linked triple algorithm \cite{How95}. However, this efficiency can only be realized in this restricted setting. 

\textbf{UACalc} also implements Freese's algorithm, and was done under his supervision in Java and Jython. It includes a GUI version in Java and a command line interface in Jython. The performance of the GUI version is poor and it will not be used for comparison. On the other hand the command line interface has a performance that is comparable to \textbf{CreamAllCongruences}. The same random Rees Matrix Semigroups that were used with \textbf{CongruencesOfSemigroup} and \textbf{CreamAllCongruences} were written into files readable by \textbf{UACalc} and the congruences were calculated in \textbf{UACalc}. \textbf{CreamAllCongruences} is consistently more than 3 times faster than \textbf{UACalc}. For these algebras with one binary operation the runtime rises roughly proportionally to $n^4$ (where $n$ is the size of the algebra) or $t^2$ (being $t$ the number cells of the multiplication matrix of the algebra operation). 

To further compare the performance of \textbf{CreamAllCongruences} with \textbf{UACalc}, tests with several specific types of Semigroups, Monoids and Groups were run. 
\begin{table}[H]
\caption{Performance Comparison between CREAM and UACalc}\label{reesSemigroupPerformance}
	\renewcommand\arraystretch{1.5}
	\noindent\[
{\tiny
	\begin{tabular}{|c|c|c|c|c|}
	\hline
    Algebra&Size&Number of&CreamAllCongruences&UACalc\\
    &&Congruences&(ms)&(ms)\\
    \hhline{|=|=|=|=|=|}
    GossipMonoid(3)&11&84&1&2\\
    \hline
    PlanarPartitionMonoid(2)&14&9&1&1\\
    \hline
    JonesMonoid(4)&14&9&$<1$&1\\
    \hline
    BrauerMonoid(3)&15&7&$<1$&2\\
    \hline
    PartitionMonoid(2)&15&13&$<1$&1\\
    \hline
    FullPBRMonoid(1)&16&167&3&4\\
    \hline
    SymmetricGroup(4)&24&4&3&10\\
    \hline
    FullTransformationSemigroup(3)&27&7&3&14\\
    \hline
    FullTransformationMonoid(3)&27&7&3&7\\
    \hline
    SymmetricInverseMonoid(3)&34&7&6&25\\
    \hline
    JonesMonoid(5)&42&6&12&36\\
    \hline
    MotzkinMonoid(3)&51&10&23&73\\
    \hline
    PartialTransformationMonoid(3)&64&7&63&223\\
    \hline
    PartialBrauerMonoid(3)&76&16&117&351\\
    \hline
    BrauerMonoid(4)&105&19&430&1 420\\
    \hline
    SymmetricGroup(5)&120&3&812&2 700\\
    \hline
    PlanarPartitionMonoid(3)&132&10&1 094&3 490\\
    \hline
    JonesMonoid(6)&132&10&950&3 444\\
    \hline
    PartitionMonoid(3)&203&16&6 500&19 979\\
    \hline
    SymmetricInverseMonoid(4)&209&11&6 970&22 318\\
    \hline
    FullTransformationSemigroup(4)&256&11&22 529&53 557\\
    \hline
    FullTransformationMonoid(4)&256&11&22 385&52 361\\
    \hline
    MotzkinMonoid(4)&323&11&47 593&134 945\\
    \hline
    JonesMonoid(7)&429&7&161 987&485 004\\
	\hline
	\end{tabular}
}
	\]
\end{table}

The performance for these algebras is mostly similar to what was found for the Rees Matrix Semigroups, with a bigger variance in the results. On average \textbf{CreamAllCongruences} is between 3 and 3,5 times faster than \textbf{UACalc}  but, depending on the algebras, we get improvements spanning from 2 times to 4,5 times. The runtime also rises roughly proportionally to $n^4$ (although with more variance).

\subsection{Automorphisms algorithm performance}
\label{subsection:automorphismsAlgPerfAnnex}

The CREAM package automorphism function was run on algebras of Type ($2^m$,$1^n$) to compare the timings against the Loops package.  
All experiments were run 7 times, with the lowest and highest times discarded. The averages of the remaining 5 are reported in Table \ref{magmautvsLoopsAnnex}.  
We see that the speeds from both packages are quite close. 

\begin{table}[H]
	\caption{Performance Comparison between Loops and CREAM on Automorphism Group Generation.}\label{magmautvsLoopsAnnex}
	\renewcommand\arraystretch{1.5}
	\noindent\[
{\tiny	\begin{tabular}{|c|c|c|}
	\hline
	Algebraic Structure&Loops Package(sec)&CREAM Package(sec)\\
	\hhline{|=|=|=|}
	Quasigroups, order 5 & 0.588 & 0.562\\
	\hline
	Quasigroups, order 6 & 520 & 541\\
	\hhline{|=|=|=|}
	Loops, order 6 & 0.101 & 0.086\\
	\hline
	Loops, order 7 & 14.623 & 13.945\\
	\hhline{|=|=|=|}
	Groups, order 32 & 0.854 & 0.935\\
	\hline
	Groups, order 64 & 39.98 & 41.05\\
	\hline
	Groups, order 128 & 12,031 & 12,160\\
	\hline
	\end{tabular}
}	\]
\end{table}

We have also run experiments on generating automorphism groups for magmas and semigroups.  Results are displayed on Table \ref{magmasemi}. For these experiments, we do not have results from other GAP packages to compare with. Therefore the following results can be taken as the benchmark to beat in future improvements.  
\begin{table}[h]
	\caption{Performance of CREAM on Automorphism Group Generation for Magmas and Semigroups.}\label{magmasemi}
	\renewcommand\arraystretch{1.5}
	\noindent\[
{\tiny	\begin{tabular}{|c|c|}
	\hline
	Algebraic Structure&Time(sec)\\
	\hhline{|=|=|}
	Magmas, order 2 & 0.016\\
	\hline
	Magmas, order 3 & 1.253\\
	\hhline{|=|=|}
	Semigroups, order 6 & 2.316\\
	\hline
	Semigroups, order 7 & 45.947\\
	\hline
	\end{tabular}
}
	\]
\end{table}

%
%
%
%
%
%
%

The CREAMAutomorphisms runs substantially faster than the simple intersection of automorphism groups, especially when the automorphism group is the trivial group. Sample results are shown in Table~\ref{table:algebra_Auto}. 

\begin{table}[H]
	\caption{Performance Comparison between CREAMAutomorphisms and Intersection of Automorphism Groups.}\label{table:algebra_Auto}
	\renewcommand\arraystretch{1.5}
	\noindent\[
{\tiny
	\begin{tabular}{|c|r|r|c|}
		\hline
		Algebra of type (2, 2, 2, ...)&\makecell{\\ Algebra\\\ AutomorphismGroup \\ (sec)}&\makecell{Intersection of \\ Automorphism \\ Groups\\(sec)}&\makecell{\\ Is Trival Group \\ the Automorphism \\Group?}\\
		\hhline{|=|=|=|=|}
		\makecell{All 2,328 non-isomorphic \\groups of order 128} & 153.8 & 201.5&No\\
		\hline
		\makecell{All 67 non-isomorphic \\groups of order 243} & 10.3 & 56.3&No\\
		\hline
		\makecell{All 15 non-isomorphic \\groups of order 625}& 11.3 & 122.9&No\\
		\hline
		\makecell{All 15 non-isomorphic \\groups of order 999} & 2.1 & 168.9&Yes\\
		\hline
	\end{tabular}
}
	\]
\end{table}

\clearpage
\subsection{Endomorphisms algorithm performance}
\label{subsection:endomorphismsAlgPerfAnnex}

To evaluate the performance of the algorithms to calculate endomorphisms,  we  tested both the classic algorithm and the congruence algorithm running with \textbf{Mace4}, as described in \ref{subsection:endomorphismsAlg}, on several types of Semigroups, Monoids and Groups. The only GAP function to calculate endomorphisms that came to our attention is the function Endomorphisms from the package SONATA. This function can calculate endomorphisms for a very narrow set of algebras, namely groups and near-rings. The groups in the list algebras were also run with SONATA.\\ 

\begin{table}[H]
	\caption{Performance Comparison between endomorphisms calculation algorithms}\label{endoPerformanceAnnex}
	\renewcommand\arraystretch{1.5}
	\noindent\[
{\tiny
	\begin{tabular}{|c|c|c|c|c|c|}
	\hline
    Algebra&Size&Number of&Classic&Congruences&SONATA\\
    &&Endomorphisms&(ms)&(ms)&(ms)\\
    \hhline{|=|=|=|=|=|=|}
    GossipMonoid(3)&11&66&30&1 635&NA\\
    \hline
    PlanarPartitionMonoid(2)&14&72&56&233&NA\\
    \hline
    JonesMonoid(4)&14&72&45&232&NA\\
    \hline
    BrauerMonoid(3)&15&28&38&160&NA\\
    \hline
    PartitionMonoid(2)&15&89&50&313&NA\\
    \hline
    FullPBRMonoid(1)&16&1 426&585&5 134&NA\\
    \hline
    SymmetricGroup(4)&24&58&101&142&166\\
    \hline
    FullTransformationSemigroup(3)&27&40&116&364&NA\\
    \hline
    FullTransformationMonoid(3)&27&40&112&360&NA\\
    \hline
    SymmetricInverseMonoid(3)&34&54&274&504&NA\\
    \hline
    JonesMonoid(5)&42&113&746&777&NA\\
    \hline
    MotzkinMonoid(3)&51&98&1 972&2 400&NA\\
    \hline
    PartialTransformationMonoid(3)&64&138&2 586&1 963&NA\\
    \hline
    PartialBrauerMonoid(3)&76&165&7 897&7 796&NA\\
    \hline
    BrauerMonoid(4)&105&274&50 164&19 875&NA\\
    \hline
    SymmetricGroup(5)&120&146&32 861&2 979&201\\
    \hline
    PlanarPartitionMonoid(3)&132&393&95 460&25 803&NA\\
    \hline
    JonesMonoid(6)&132&393&131 472&22 889&NA\\
    \hline
    PartitionMonoid(3)&203&687&741 441&105 421&NA\\
    \hline
    SymmetricInverseMonoid(4)&209&282&470 901&75 940&NA\\
	\hline
	\end{tabular}
}
	\]
\end{table}

In this test the classic algorithm not using congruences is faster than using congruences for algebras up to size 60 but for larger algebras its runtime raises very fast and the algorithm using congruences becomes faster.\\

In view of these results, the \textbf{CreamEndomorphisms} function uses the classic algorithm for algebras with size up to 60 and the congruences algorithm for larger algebras.\\


\clearpage
\section{Applications}
\label{section:applications}

\subsection{Monolithic Algebras}
\label{subsection:monolithic}

Some examples of the use of the package to calculate whether an algebra is monolithic can be seen next:

{\tiny
\begin{verbatim}
gap> algebra := [RecoverMultiplicationTable(6,1)];
[ [ [ 1, 1, 1, 1, 1, 1 ], [ 1, 1, 1, 1, 1, 1 ], [ 1, 1, 1, 1, 1, 1 ], 
      [ 1, 1, 1, 1, 1, 1 ], [ 1, 1, 1, 1, 1, 1 ], [ 1, 1, 1, 1, 1, 1 ] ] ]
gap> CreamAllPrincipalCongruences(algebra);
[ [ -1, -1, -1, -1, -2, 5 ], [ -1, -1, -1, -2, 4, -1 ], 
  [ -1, -1, -1, -2, -1, 4 ], [ -1, -1, -2, 3, -1, -1 ], 
  [ -1, -1, -2, -1, 3, -1 ], [ -1, -1, -2, -1, -1, 3 ], 
  [ -1, -2, 2, -1, -1, -1 ], [ -1, -2, -1, 2, -1, -1 ], 
  [ -1, -2, -1, -1, 2, -1 ], [ -1, -2, -1, -1, -1, 2 ], 
  [ -2, 1, -1, -1, -1, -1 ], [ -2, -1, 1, -1, -1, -1 ], 
  [ -2, -1, -1, 1, -1, -1 ], [ -2, -1, -1, -1, 1, -1 ], 
  [ -2, -1, -1, -1, -1, 1 ] ]
gap> CreamIsAlgebraMonolithic(algebra);
false
gap> algebra := [RecoverMultiplicationTable(6,19)];
[ [ [ 1, 1, 1, 1, 1, 1 ], [ 1, 1, 1, 1, 1, 1 ], [ 1, 1, 1, 1, 1, 1 ], 
      [ 1, 1, 1, 1, 1, 1 ], [ 1, 1, 1, 2, 1, 1 ], [ 1, 1, 2, 1, 1, 1 ] ] ]
gap> CreamAllPrincipalCongruences(algebra);        
[ [ -2, 1, -1, -1, -1, -1 ], [ -2, 1, -1, -1, -2, 5 ], 
  [ -2, 1, -1, -2, 4, -1 ], [ -2, 1, -1, -2, -1, 4 ], 
  [ -2, 1, -2, 3, -1, -1 ], [ -2, 1, -2, -1, 3, -1 ], 
  [ -2, 1, -2, -1, -1, 3 ], [ -3, 1, 1, -1, -1, -1 ], 
  [ -3, 1, -1, 1, -1, -1 ], [ -3, 1, -1, -1, 1, -1 ], 
  [ -3, 1, -1, -1, -1, 1 ] ]
gap> CreamIsAlgebraMonolithic(algebra);            
true
\end{verbatim}
}

\subsection{Small Semigroups}
\label{subsection:smallsemi}

For all small semigroups up to size 6, all its congruences and endomorphisms were calculated. 
It was moreover determined whether the semigroups were monolithic (this was calculated up to size 8).

\begin{table}[H]
	\caption{Determination of all monolithic semigroups up to size 8 }\label{monolithic performance}
	\renewcommand\arraystretch{1.5}
	\noindent\[
{\tiny
	\begin{tabular}{|c|c|c|c|}
	\hline
    \textbf{Size}&\textbf{Number of}&\textbf{Number of}&\textbf{Performance}\\
    &\textbf{semigroups}&\textbf{monoliths}&\textbf{(ms)}\\
    \hline
    \hline
    1&1&0&0\\
    \hline
    2&4&4&0\\
    \hline
    3&18&7&0\\
    \hline
    4&126&16&8\\
    \hline
    5&1 160&103&108\\
    \hline
    6&15 973&1 823&1 712\\
    \hline
    7&836 021&149 020&127 356\\
    \hline
    8&1 843 120 128&48 438 046&462 897 348\\
    \hline
	\end{tabular}
}
	\]
\end{table}

All automorphisms for semigroups up to size 7 were also calculated separately and the automorphism group and its ID was calculated for each of these semigroups. The group ID is the identification of a group in the Small Group GAP library \cite{SMALLGRP1_4_2}.

The CREAM library doesn't calculate the automorphisms of an algebra as an automorphism group but as a list of bijective mappings, but there is an easy way to calculate the automorphism group from the list of mappings. To do this, all the mappings are converted into permutations using the \textbf{PermList} function and each of these permutations is used as generator for the group. Having this we can get the Group Id using the \textbf{IdGroup} function. 
%

\clearpage
The following table gives the automorphism of semigroups up to size $7$ (and reproduces the results in \cite{joao3}):
{\tiny
\begin{table}[h]
	\caption{Computation of the automorphism groups of semigroups up to size 7 }\label{automorphism smallsemi}
	\renewcommand\arraystretch{1.5}
	\noindent\[
{\tiny
	\begin{tabular}{|l|c|c|c|c|c|c|}
	\hline
    \textbf{Size}&2&3&4&5&6&7\\
    \hline
    \textbf{\# Semigroups}&4&18&126&1 160&15 973&836 021\\
    \hline
    \hline
    [ 1, 1 ] trivial&3&12&78&746&10 965&746 277\\
    \hline
    [ 2, 1 ] $C_2$&1&5&39&342&4 121&76 704\\
    \hline
    [ 3, 1 ] $C_3$&&&&2&26&412\\
    \hline
    [ 4, 1 ] $C_4$&&&&1&7&82\\
    \hline
    [ 4, 2 ] $C_2 \times C_2$&&&3&26&441&7 314\\
    \hline
    [ 5, 1 ] $C_5$&&&&&&6\\
    \hline
    [ 6, 1 ] $S_3$&&1&5&33&300&3 638\\
    \hline
    [ 6, 2 ] $C_6$&&&&&&37\\
    \hline
    [ 8, 2 ] $C_4 \times C_2$&&&&&&4\\
    \hline
    [ 8, 3 ] $D_8$&&&&1&17&169\\
    \hline
    [ 8, 5 ] $C_2 \times C_2 \times C_2$&&&&&6&172\\
    \hline
    [ 10, 1 ] $D_{10}$&&&&&&2\\
    \hline
    [ 12, 4 ] $D_{12}$&&&&4&49&790\\
    \hline
    [ 16, 11 ] $C_2 \times D_8$&&&&&&10\\
    \hline
    [ 24, 12 ] $S_4$&&&1&4&30&277\\
    \hline
    [ 24, 14 ] $C_2 \times C_2 \times S_3$&&&&&&14\\
    \hline
    [ 36, 10 ] $S_3 \times S_3$&&&&&2&24\\
    \hline
    [ 48, 48 ] $C_2 \times S_4$&&&&&4&45\\
    \hline
    [ 72, 40 ] $(S_3 \times S_3) \wr C_2$&&&&&&1\\
    \hline
    [ 120, 34 ] $S_5$&&&&1&4&30\\
    \hline
    [ 144, 183 ] $S_3 \times S_4$&&&&&&4\\
    \hline
    [ 240, 189 ] $C_2 \times S_5$&&&&&&4\\
    \hline
    [ 720, 763 ] $S_6$&&&&&1&4\\
    \hline
    [ 5040, - ] $S_7$&&&&&&1\\
    \hline
	\end{tabular}
}
	\]
\end{table}
 }

\clearpage
\subsection{Small Groups}
\label{subsection:smallgroups}

Using the Small Group GAP library \cite{SMALLGRP1_4_2} the groups from order 2 to 96 were extracted and the Cream functions were used with them:

\begin{itemize}
    \item Order 2-31 - CreamAllCongruences, CreamAllEndomorphisms and CreamIsMonolithic
    \item Order 32-96 - CreamAllCongruences and CreamIsMonolithic
\end{itemize}


A summary from these runs is provided in Table~\ref{smallgroupruns}.\\
\begin{table}[h]
	\caption{Summary for Small Group Runs}\label{smallgroupruns}
	\renewcommand\arraystretch{1.5}
	\noindent\[
{\tiny
	\begin{tabular}{|c|c|c|}
	\hline
    \textbf{Size}&\textbf{Number of}&\textbf{Number of}\\
    &\textbf{Groups}&\textbf{Monolithic}\\
    \hline
    \hline
    2-7&8&6\\
    \hline
    8-15&19&9\\
    \hline
    16&14&6\\
    \hline
    17-23&17&7\\
    \hline
    24&15&2\\
    \hline
    25-31&19&7\\
    \hline
    32&50&16\\
    \hline    
    33-47&54&10\\
    \hline
    48&52&4\\
    \hline
    49-63&69&16\\
    \hline
    64&52&22\\
    \hline
    65-71&17&3\\
    \hline
    72&50&3\\
    \hline
    73-79&18&5\\
    \hline
    80&52&1\\
    \hline
    81-95&70&13\\
    \hline
    96&231&13\\
    \hline
	\end{tabular}
}
	\]
\end{table}

\clearpage
\subsection{Groups, Semigroups and Monoids}
\label{subsection:largesemi}

For a selection of larger groups, semigroups and monoids, all congruences and endomorphisms were calculated. It was also calculated whether these algebras were monolithic. A summary from these runs is provided in Table~\ref{selectalgebras}.

\begin{table}[h]
	\caption{\# of Congruences and Endomorphisms for a selection of larger algebras }\label{selectalgebras}
	\renewcommand\arraystretch{1.5}
	\noindent\[
{\tiny
	\begin{tabular}{|c|c|c|c|c|}
	\hline
    \textbf{Algebra}&\textbf{Size}&\textbf{Number of}&\textbf{Number of}&\textbf{Is}\\
    &&\textbf{Congruences}&\textbf{Endomorphisms}&\textbf{Monolithic}\\
    \hline
    \hline
    GossipMonoid(3)&11&84&66&No\\
    \hline
    PlanarPartitionMonoid(2)&14&9&72&No\\
    \hline
    JonesMonoid(4)&14&9&72&No\\
    \hline
    BrauerMonoid(3)&15&7&28&No\\
    \hline
    PartitionMonoid(2)&15&13&89&No\\
    \hline
    FullPBRMonoid(1)&16&167&1426&No\\
    \hline
    SymmetricGroup(4)&24&4&58&Yes\\
    \hline
    FullTransformationSemigroup(3)&27&7&40&Yes\\
    \hline
    FullTransformationMonoid(3)&27&7&40&Yes\\
    \hline
    SymmetricInverseMonoid(3)&34&7&54&Yes\\
    \hline
    JonesMonoid(5)&42&6&113&No\\
    \hline
    MotzkinMonoid(3)&51&10&98&No\\
    \hline
    PartialTransformationMonoid(3)&64&7&138&Yes\\
    \hline
    PartialBrauerMonoid(3)&76&16&165&No\\
    \hline
    BrauerMonoid(4)&105&19&274&No\\
    \hline
    SymmetricGroup(5)&120&3&146&Yes\\
    \hline
    PlanarPartitionMonoid(3)&132&10&393&No\\
    \hline
    JonesMonoid(6)&132&10&393&No\\
    \hline
    PartitionMonoid(3)&203&16&687&No\\
    \hline
    SymmetricInverseMonoid(4)&209&11&282&Yes\\
    \hline
    FullTransformationSemigroup(4)&256&11&345&Yes\\
    \hline
    FullTransformationMonoid(4)&256&11&345&Yes\\
    \hline
	\end{tabular}
}
	\]
\end{table}

\clearpage
\subsection{Automorphisms of Specific Algebra Classes}
\label{subsection:autoalgebraclasses}

 The following table contains the groups that appear as automorphims groups of associative rings up to order $15$.

\begin{table}[h]
\caption{Computation of the automorphism groups of rings up to size 15 }\label{automorphism smallring}
\renewcommand\arraystretch{1.5}
\noindent\[
{\tiny
	\begin{tabular}{|l|c|c|c|c|c|c|c|c|c|c|c|c|c|c|}
	\hline
	\textbf{Size}&2&3&4&5&6&7&8&9&10&11&12&13&14&15\\
	\hline
	\textbf{\# Rings}&2&2&11&2&4&2&52&11&4&2&22&2&4&4\\
	\hline
	\hline
	[ 1, 1 ] trivial&2&1&3& &2& &4& & & &3& & &  \\
	\hline
	[ 2, 1 ] $C_2$& &1&6&1& &1&12&1&2& &4& &2&1\\
	\hline
	[ 4, 1 ] $C_4$&&&&&&&&&&1&&&& \\
	\hline
	[ 4, 2 ] $C_2 \times C_2$&&&&&1&&6&3&&&&1&&\\
	\hline
	[ 6, 1 ] $S_3$&&&2&&&&1&&&&&&&\\
	\hline
	[ 6, 2 ] $C_6$&&&&&&&1&&&&&&&\\
	\hline
	[ 8, 3 ] $D_8$&&&&&&&3&&&&&&&\\
	\hline
	[ 8, 5 ] $C_2 \times C_2 \times C_2$&&&&&&&4&1&&&2&&&\\
	\hline
	[ 12, 4 ] $D_{12}$&&&&&&&&2&&&&&&\\
	\hline
	[ 16, 11 ] $C_2 \times D_8$&&&&&&&3&&&&&&&\\
	\hline
	[ 16, 14 ] $C_2 \times C_2 \times C_2 \times C_2$&&&&&&&&&&&2&&&\\
	\hline
	[ 24, 12 ] $S_4$&&&&&1&&&&&&&&&\\
	\hline
	[ 24, 13 ] $C_2 \times A_4 \times S_3$&&&&&&&2&&&&&&&\\
	\hline
	[ 32, 46 ] $C_2 \times C_2 \times D_8$&&&&&&&&&&&4&&&\\
	\hline
	[ 36, 10 ] $S_3 \times S_3$&&&&&&&4&&&&&&&\\
	\hline
	[ 48, 48 ] $C_2 \times S_4$&&&&&&&7&&&&&&&\\
	\hline
	[ 64, 261 ] $C_2 \times C_2 \times C_2 \times D_8$&&&&&&&&&&&&&&1\\
	\hline
	[ 72, 40 ] $(S_3 \times S_3) \wr C_2$&&&&&&&&2&&&&&&\\
	\hline
	[ 120, 34 ] $S_5$&&&&&1&&&&&&&&&\\
	\hline
	[ 144, 183 ] $S_3 \times S_4$&&&&&&&2&&&&&&&\\
	\hline
	[ 216, 162 ] $S_3 \times S_3 \times S_3$&&&&&&&&&&&2&&&\\
	\hline
	[ 576, 8653 ] $S_4 \times S_4$&&&&&&&&&1&&&&&\\
	\hline
	[ 720, 763 ] $S_6$&&&&&&1&&&&&&&&\\
	\hline
	[ 5040, - ] $S_7$&&&&&&&3&&&&&&&\\
	\hline
	[ 13824, - ] $S_4 \times S_4\times S_4$&&&&&&&&&&&&&&1\\
	\hline
	[ 14400, - ] $S_5 \times S_5$&&&&&&&&&&&1&&&\\
	\hline
	[ 17280, - ] $S_6 \times S_4$&&&&&&&&&&&2&&&\\
	\hline
	[ 40320, - ] $S_8$&&&&&&&&2&&&&&&\\
	\hline
	[ 362880, - ] $S_{10}$&&&&&&&&&1&1&&&&\\
	\hline
	[ 518400, - ] $S_6 \times S_6$&&&&&&&&&&&&&1&\\
	\hline
	[ 39916800, - ] $S_{11}$&&&&&&&&&&&2&&&\\
	\hline
	[ 479001600, - ] $S_{12}$&&&&&&&&&&&&1&&\\
	\hline
	[ 6227020800, - ] $S_{13}$&&&&&&&&&&&&&1&\\
	\hline
	[ 87178291200, - ] $S_{14}$&&&&&&&&&&&&&&1\\
	\hline
	\end{tabular}
}
\]
\end{table}

\clearpage
The variety of {\em Quasi-MV-algebras} \cite{quasiMV} is defined by the following identities (in ProverX syntax \cite{araujo19c,proverx}): 
{\tiny
\begin{verbatim}
(x + z) + y = x + (y + z).
x'' = x.
x + 1 = 1.
(x' + y)' + y = (y' + x)' + x.
(x + 0)' = x' + 0.
(x + y) + 0 = x + y.
0' = 1.
\end{verbatim}  
}
\begin{table}[h]
	\caption{Computation of the automorphism groups of Quasi-MV-algebras of size up to 12 }\label{automorphism small_quasiMValgebra}
	\renewcommand\arraystretch{1.5}
	\noindent\[
	{\tiny
		\begin{tabular}{|l|c|c|c|c|c|c|c|c|c|c|c|}
		\hline
		\textbf{Size}&2&3&4&5&6&7&8&9&10&11&12\\
		\hline
		\textbf{\# Quasi-MV-algebras}&1&1&4&4&11&11&27&27&60&62&131\\
		\hline
		\hline
		[ 1, 1 ] trivial&1&1&3&2&7&4&16&8&35&17&76\\
		\hline
		[ 2, 1 ] $C_2$&&&1&2&3&4&5&9&10&18&19\\
		\hline
		[ 4, 2 ] $C_2 \times C_2$&&&&&&2&1&4&2&10&5\\
		\hline
		[ 6, 1 ] $S_3$&&&&&1&&3&&5&&9\\
		\hline
		[ 8, 5 ] $C_2 \times C_2 \times C_2$&&&&&&&&2&1&4&2 \\
		\hline
		[ 12, 4 ] $D_{12}$&&&&&&&1&&2&&5\\
		\hline
		[ 16, 14 ] $C_2 \times C_2 \times C_2 \times C_2$&&&&&&&&&&2&1 \\
		\hline
		[ 24, 12 ] $S_4$&&&&&&1&&2&&4&\\
		\hline
		[ 24, 14 ] $C_2 \times C_2 \times C_3$&&&&&&&&&1&&2\\
		\hline
		[ 48, 48 ] $C_2 \times S_4$&&&&&&&&1&&2&\\
		\hline
		[ 48, 51 ] $C_2 \times C_2 \times C_2 \times S_3$&&&&&&&&&&&1\\
		\hline
		[ 96, 226 ] $C_2 \times C_2 \times S_4$&&&&&&&&&&1&\\
		\hline
		[ 120, 34 ] $S_5$&&&&&&&1&1&2&&4\\
		\hline
		[ 240, 189 ] $C_2 \times S_5$&&&&&&&&&1&&2\\
		\hline
		[ 480, 1186 ] $C_2 \times C_2 \times S_5$&&&&&&&&&&&1\\
		\hline
		[ 720, 763 ] $S_6$&&&&&&&&&&2&\\
		\hline
		[ 1440, 5842 ] $C_2 \times S_6$&&&&&&&&&&1&\\
		\hline
		[ 5040, - ] $S_7$&&&&&&&&&1&&2\\
		\hline
		[ 10080, - ] $C_2 \times S_7$&&&&&&&&&&&1\\
		\hline
		[ 40320, - ] $S_8$&&&&&&&&&&1&\\
		\hline
		[ 362880, - ] $S_9$&&&&&&&&&&&1\\
		\hline
		\end{tabular}
	}
	\]
\end{table}

\clearpage
The quasivariety of {\em BCI algebras} \cite{bci} is defined by the following identities (in ProverX syntax \cite{araujo19c,proverx}):
{\tiny
\begin{verbatim}
((x * y) * (x * z)) * (z * y) = 0.
(x * (x * y)) * y = 0.
x * x = 0.
((x * y = 0) & (y * x = 0)) -> (x = y).
(x * 0 = 0) -> (x = 0).
\end{verbatim}
}

\begin{table}[h]
	\caption{Computation of the automorphism groups of BCI algebras of size up to 7 }\label{automorphism small_BCIalgebra}
	\renewcommand\arraystretch{1.5}
	\noindent\[
	{\tiny
		\begin{tabular}{|l|c|c|c|c|c|c|}
		\hline
		\textbf{Size}&2&3&4&5&6&7\\
		\hline
		\textbf{\# BCI algebras}&2&5&22&118&974&10,834\\
		\hline
		\hline
		[ 1, 1 ] trivial&2&3&13&78&679&7,970\\
		\hline
		[ 2, 1 ] $C_2$&&2&7&31&241&2,384\\
		\hline
		[ 4, 1 ] $C_4$&&&&1&1&2\\
		\hline
		[ 4, 2 ] $C_2 \times C_2$&&&&2&20&207\\
		\hline
		[ 6, 1 ] $S_3$&&&2&5&25&195\\
		\hline
		[ 6, 2 ] $C_6$&&&&&&1\\
		\hline
		[ 8, 2 ] $C_4 \times C_2$&&&&&&1\\
		\hline
		[ 8, 3 ] $D_8$&&&&&&4\\
		\hline
		[ 8, 5 ] $C_2 \times C_2 \times C_2$&&&&&&3\\
		\hline
		[ 12, 4 ] $D_{12}$&&&&&4&35\\
		\hline
		[ 24, 12 ] $S_4$&&&&1&4&22\\
		\hline
		[ 36, 10 ] $S_3 \times S_3$&&&&&&2\\
		\hline
		[ 48, 48 ] $C_2 \times S_4$&&&&&&3\\
		\hline
		[ 120, 34 ] $S_5$&&&&&1&4\\
		\hline
		[ 720, 763 ] $S_6$&&&&&&1\\
		\hline
		\end{tabular}
	}
	\]
\end{table}

\clearpage
The variety of {\em quandles} \cite{quandles} is defined by the following identities (in ProverX syntax \cite{araujo19c,proverx}): 
{\tiny \begin{verbatim}
x v (y v z) = (x v y) v (x v z).
(x ^ y) ^ z = (x ^ z) ^ (y ^ z).
(x v y) ^ x = y.
x v (y ^ x) = y.
x v x = x.\end{verbatim}  }

\begin{table}[h]
\caption{Computation of the automorphism groups of  quandles of size up to 6 }\label{automorphism small_Quandles}
\renewcommand\arraystretch{1.5}
\noindent\[
{\tiny
	\begin{tabular}{|l|c|c|c|c|c|}
	\hline
	\textbf{Size}&2&3&4&5&6\\
	\hline
	\textbf{\# quandles}&1&3&7&22&73\\
	\hline
	\hline
	[ 2, 1 ] $C_2$&1&1&1&&\\
	\hline
	[ 3, 1 ] $C_3$&&&1&1&\\
	\hline
	[ 4, 1 ] $C_4$&&&&1&2\\
	\hline
	[ 4, 2 ] $C_2 \times C_2$&&&1&5&8\\
	\hline
	[ 5, 1 ] $C_5$&&&&&1\\
	\hline
	[ 6, 1 ] $S_3$&&2&1&2&3\\
	\hline
	[ 6, 2 ] $C_6$&&&&2&15\\
	\hline
	[ 8, 2 ] $C_4 \times C_2$&&&&&2\\
	\hline
	[ 8, 3 ] $D_8$&&&1&3&8\\
	\hline
	[ 8, 5 ] $C_2 \times C_2 \times C_2$&&&&&7\\
	\hline
	[ 12, 3 ] $A_{4}$&&&1&1&\\
	\hline
	[ 12, 4 ] $D_{12}$&&&&3&6\\
	\hline
	[ 16, 11 ] $C_2 \times D_8$&&&&&5\\
	\hline
	[ 18, 3 ] $C_3 \times S_3$&&&&&3\\
	\hline
	[ 20, 3 ] $C_5 \times C_4$&&&&3&3\\
	\hline
	[ 24, 12 ] $S_4$&&&1&&2\\
	\hline
	[ 24, 13 ] $C_2 \times A_4 \times S_3$&&&&&3\\
	\hline
	[ 36, 10 ] $S_3 \times S_3$&&&&&1\\
	\hline
	[ 48, 48 ] $C_2 \times S_4$&&&&&2\\
	\hline
	[ 72, 40 ] $(S_3 \times S_3) \wr C_2$&&&&&1\\
	\hline
	[ 120, 34 ] $S_5$&&&&1&\\
	\hline
	[ 720, 763 ] $S_6$&&&&&1\\
	\hline
	\end{tabular}
}
\]
\end{table}

\clearpage
\subsection{Algebra Classes}
\label{subsection:algebraclasses}

The CREAM package allows for the generation of lists of algebras of specific types based on its axiomatic definition. This is done using the function \textbf{CreamAlgebrasFromAxioms} that uses Mace4 to generate algebras according with this axiomatic definition. This function takes as input a string with the axioms in Mace4 format and a positive integer with the size of the generated algebras.  The generated algebras can then be used with the CREAM functions that calculate congruences, endomorphisms and whether the algebra is monolithic.
\begin{samepage}
CREAM was applied in this way to the following classes of algebras:
\begin{itemize}
    \item Almost distributive lattices
    \item BCI-algebras
    \item Bilattices
    \item Boolean algebras
    \item Boolean rings
    \item Commutative lattice-ordered monoids
    \item Commutative regular rings
    \item Complemented distributive lattices
    \item Complemented modular lattices
    \item Distributive lattice ordered semigroups
    \item Ockham algebras
    \item Ortholattices
    \item Orthomodular lattices
    \item Idempotent semirings
    \item Extra loops
    \item Digroups
    \item Commutative dimonoids
\end{itemize}
\end{samepage}
A summary from these runs is provided in Table~\ref{algebraclass}.\\

\begin{table}[htb]
	\caption{Summary for Algebra Class Runs}\label{algebraclass}
	\renewcommand\arraystretch{1.5}
	\noindent\[
{\tiny
	\begin{tabular}{|c|c|c|c|c|}
	\hline
    \textbf{Algebra}&\textbf{Size}&\textbf{Algebra}&\textbf{Number of}&\textbf{Number of}\\
    \textbf{Class}&&\textbf{Type}&\textbf{Algebras}&\textbf{Monolithic}\\
    \hline
    \hline
    Almost distributive lattices&4&$(2^2)$&2&0\\
    \hline
    BCI-algebras&5&$(2^1)$&118&80\\
    \hline
    Bilattice&6&$(2^4,1^1)$&32&32\\
    \hline
    Boolean algebra&8&$(2^2,1^1)$&1&0\\
    \hline
    Boolean algebra&16&$(2^2,1^1)$&1&0\\
    \hline
    Boolean ring&16&$(2^2)$&1&0\\
    \hline
    Commutative lattice-ordered monoids&5&$(2^3)$&199&97\\
    \hline    
    Commutative regular rings&6&$(2^2,1^2)$&72&66\\
    \hline
    Complemented distributive lattices&16&$(2^2,1^1)$&1&0\\
    \hline
    Complemented distributive lattices&32&$(2^2,1^1)$&1&0\\
    \hline
    Complemented modular lattices&8&$(2^2,1^1)$&41&40\\
    \hline
    Distributive lattice ordered semigroups&4&$(2^3)$&479&170\\
    \hline
    Ockham algebras&6&$(2^2,1^1)$&197&20\\
    \hline
    Ortholattices&7&$(2^2,1^1)$&46&12\\
    \hline
    Orthomodular lattices&14&$(2^2,1^1)$&33&31\\
    \hline
    Idempotent semiring&5&$(2^2)$&149&42\\
    \hline
    Extra loop&10&$(2^3)$&2&1\\
    \hline
    Digroup&8&$(2^2,1^1)$&10&3\\
    \hline
    Commutative dimonoid&4&$(2^2)$&101&37\\
    \hline
	\end{tabular}
}
	\]
\end{table}

\end{appendices}

\end{document}